\input amstex.tex
\documentstyle{amsppt}
\magnification=\magstep1
\vsize 8.5truein
\hsize 6truein

\define\flow{\left(\bold{M},\{S^t\}_{t\in\Bbb R},\mu\right)}

\define\traj{S^{[a,b]}x_0}

\heading
Conditional Proof of the Boltzmann-Sinai Ergodic Hypothesis
\endheading
 
\bigskip \bigskip
 
\centerline{{\bf N\'andor Sim\'anyi}
\footnote{Research supported by the National Science Foundation, grants
DMS-0457168, DMS-0800538.}}

\bigskip \bigskip

\centerline{University of Alabama at Birmingham}
\centerline{Department of Mathematics}
\centerline{Campbell Hall, Birmingham, AL 35294 U.S.A.}
\centerline{E-mail: simanyi\@math.uab.edu}

\bigskip \bigskip 

\noindent
{\it Dedicated to Yakov G. Sinai and Domokos Sz\'asz}

\bigskip \bigskip

\hbox{\centerline{\vbox{\hsize 8cm {\bf Abstract.} We consider the system of
      $N$ ($\ge2$) elastically colliding hard balls of masses
      $m_1,\dots,m_N$ and radius $r$ on the flat unit torus $\Bbb
      T^\nu$, $\nu\ge2$. We prove the so called Boltzmann-Sinai
      Ergodic Hypothesis, i. e. the full hyperbolicity and ergodicity
      of such systems for every selection $(m_1,\dots,m_N;r)$ of the
      external parameters, provided that almost every
      singular orbit is geometrically hyperbolic (sufficient),
      i. e. the so called Chernov-Sinai Ansatz is true.  The present
      proof does not use the formerly developed, rather involved
      algebraic techniques, instead it employs exclusively dynamical
      methods and tools from geometric analysis.}}}

\bigskip \bigskip

\noindent
Primary subject classification: 37D50

\medskip

\noindent
Secondary subject classification: 34D05

\bigskip \bigskip

\heading
\S1. Introduction
\endheading

\bigskip

In this paper we prove the Boltzmann--Sinai Ergodic Hypothesis under the
condition of the Chernov-Sinai Ansatz (see \S2).
In a loose form, as attributed to L. Boltzmann back in the 1880's, 
this hypothesis asserts that gases of hard balls are ergodic. In a precise 
form, which is due to Ya. G. Sinai in 1963 [Sin(1963)], it states that the gas
of $N\ge2$ identical hard balls (of "not too big" radius) on a torus
$\Bbb T^\nu$, $\nu\ge2$ (a $\nu$-dimensional box with periodic
boundary conditions), is ergodic, provided that certain necessary
reductions have been made. The latter means that one fixes the
total energy, sets the total momentum to zero, and restricts the
center of mass to a certain discrete lattice within the torus. The
assumption of a not too big radius is necessary to have the interior of the
configuration space connected.

Sinai himself pioneered rigorous mathematical studies of hard ball
gases by proving the hyperbolicity and ergodicity for the case
$N=2$ and $\nu=2$ in his seminal paper [Sin(1970)], 
where he laid down the foundations
of the modern theory of chaotic billiards. Then Chernov and Sinai
extended this result to ($N=2$, $\nu\ge 2$), as well as proved
a general theorem on ``local'' ergodicity applicable to systems of
$N>2$ balls [S-Ch(1987)]; the latter became instrumental in the
subsequent studies. The case $N>2$ is substantially more difficult
than that of $N=2$ because, while the system of two balls reduces
to a billiard with strictly convex (spherical) boundary, which
guarantees strong hyperbolicity, the gases of $N>2$ balls reduce
to billiards with convex, but not strictly convex, boundary (the
latter is a finite union of cylinders) -- and those are characterized
by very weak hyperbolicity.

Further development has been due mostly to A. Kr\'amli,
D. Sz\'asz, and the present author. We proved hyperbolicity and
ergodicity for $N=3$ balls in any dimension [K-S-Sz(1991)] by 
exploiting the ``local'' ergodic theorem of Chernov and Sinai 
[S-Ch(1987)], and carefully analyzing all
possible degeneracies in the dynamics to obtain ``global''
ergodicity. We extended our results to $N=4$ balls in
dimension $\nu\ge3$ next year [K-S-Sz(1992)], and then I proved the
ergodicity whenever $N\le\nu$ [Sim(1992)-I-II] (this covers systems with an
arbitrary number of balls, but only in spaces of high enough
dimension, which is a restrictive condition). At this point,
the existing methods could no longer handle any new cases, because
the analysis of the degeneracies became overly complicated. It was
clear that further progress should involve novel ideas.

A breakthrough was made by Sz\'asz and myself, when we used the methods of
algebraic geometry [S-Sz(1999)]. We assumed that the balls had arbitrary
masses $m_1,\dots,m_N$ (but the same radius $r$). Now by taking the limit
$m_N\to 0$, we were able to reduce the dynamics of $N$ balls to the motion of
$N-1$ balls, thus utilizing a natural induction on $N$. Then algebro-geometric
methods allowed us to effectively analyze all possible degeneracies, but only
for typical (generic) $(N+1)$-tuples of ``external'' parameters
$(m_1,\dots,m_N,r)$; the latter needed to avoid some exceptional submanifolds
of codimension one, which remained unknown. This approach led to a proof of
full hyperbolicity (but not yet ergodicity!) for all $N\ge2$ and $\nu\ge2$, and
for generic $(m_1,\dots,m_N,r)$, see [S-Sz(1999)]. Later the present author
simplified the arguments and made them more ``dynamical'', which allowed me to
obtain full hyperbolicity for hard balls with any set of external
parameters $(m_1,\dots,m_N,r)$ [Sim(2002)]. (The reason why the masses $m_i$
are considered {\it geometric parameters} is that they determine the relevant
Riemannian metric

$$
||dq||^2=\sum_{i=1}^N m_i||dq_i||^2
$$
of the system, see \S2 below.) Thus, the hyperbolicity has been fully 
established for all systems of hard balls on tori.

To upgrade the full hyperbolicity to ergodicity, one needs to refine the
analysis of the aforementioned degeneracies. For hyperbolicity, it was enough
that the degeneracies made a subset of codimension $\ge1$ in the phase space.
For ergodicity, one has to show that its codimension is $\ge2$, or to find
some other ways to prove that the (possibly) arising codimension-one manifolds
of non-sufficiency are incapable of separating distinct ergodic components.
The latter approach will be pursued in this paper. In the paper [Sim(2003)] I
took the first step in the direction of proving that the codimension of
exceptional manifolds is at least two: I proved that the systems of $N\ge2$
disks on a 2D torus (i.e., $\nu=2$) are ergodic for typical (generic)
$(N+1)$-tuples of external parameters $(m_1,\dots,m_N,r)$. The proof again
involves some algebro-geometric techniques, thus the result is restricted to
generic parameters $(m_1,\dots,m_N;\,r)$.  But there was a good reason to
believe that systems in $\nu\ge3$ dimensions would be somewhat easier to
handle, at least that was indeed the case in early studies.

Finally, in my paper [Sim(2004)] I was able to further improve the
algebro-geometric methods of [S-Sz(1999)], and proved that for any $N\ge2$,
$\nu\ge2$ and for almost every selection $(m_1,\dots,m_N;\,r)$ of the external
geometric parameters the corresponding system of $N$ hard balls on 
$\Bbb T^\nu$ is (fully hyperbolic and) ergodic.

\medskip

In this paper I will prove the following result.

\medskip

\subheading{Theorem} For any integer values $N\ge2$, $\nu\ge2$,
and for every $(N+1)$-tuple $(m_1,\dots,m_N,r)$ of the external geometric
parameters the standard hard ball system 
$\left(\bold M_{\vec m,r},\,\left\{S_{\vec m,r}^t\right\},\,
\mu_{\vec m,r}\right)$ is (fully hyperbolic and) ergodic, provided
that the Chernov-Sinai Ansatz (see \S2) holds true for all such systems.

\medskip

\subheading{Remark 1.1} The novelty of the theorem (as compared to the result 
in [Sim(2004)]) is that it applies to each $(N+1)$-tuple of external 
parameters (provided that the interior of the phase space is connected), 
without an exceptional zero-measure set. Somehow the most annoying shortcoming
of several earlier results has been exactly the fact that those results are
only valid for hard sphere systems apart from an undescribed, countable
collection of smooth, proper submanifolds of the parameter space 
$\Bbb R^{N+1}\ni(m_1,m_2,\dots,m_N;\,r)$. Furthermore, those proofs do not
provide any effective means to check if a given $(m_1,\dots,m_N;\,r)$-system
is ergodic or not, most notably for the case of equal masses in Sinai's
classical formulation of the problem.

\medskip

\subheading{Remark 1.2} The present result speaks about exactly the same 
models as the result of [Sim(2002)], but the statement of this new theorem is 
obviously stronger than that of the theorem in [Sim(2002)]: It has been known
for a long time that, for the family of semi-dispersive billiards, ergodicity
cannot be obtained without also proving full hyperbolicity.

\medskip

\subheading{Remark 1.3} As it follows from the results of [C-H(1996)] and
[O-W(1998)], all standard hard ball systems $\flow$ (the models covered by the
theorem) are not only ergodic, but they enjoy the Bernoulli mixing property, 
as long as they are known to be mixing. However, even the K-mixing property 
of semi-dispersive billiard systems follows from their ergodicity, as the
classical results of Sinai in [Sin(1968)], [Sin(1970)], and [Sin(1979)] show.

\bigskip

\subheading{The Organization of the Paper} 
In the subsequent section we overview the necessary technical
prerequisites of the proof, along with many of the needed references
to the literature. The fundamental objects of this paper are the so
called "exceptional manifolds" or "separating manifolds" $J$: they are
codimension-one submanifolds of the phase space that are separating
distinct, open ergodic components of the billiard flow.

\medskip

In \S3 we prove our Main Lemma (3.5), which states, roughly speaking,
the following: Every separating manifold $J\subset\bold M$ contains at
least one sufficient (or geometrically hyperbolic, see \S2) phase
point. The existence of such a sufficient phase point $x\in J$,
however, contradicts the Theorem on Local Ergodicity of Chernov and
Sinai (Theorem 5 in [S-Ch(1987)]), for an open neighborhood $U$ of $x$
would then belong to a single ergodic component, thus violating the
assumption that $J$ is a separating manifold. In \S4 this result will
be exploited to carry out an inductive proof of the (hyperbolic)
ergodicity of every hard ball system, provided that the Chernov-Sinai
Ansatz (see \S2) holds true for all hard ball systems.

In what follows, we make an attempt to briefly outline the key ideas
of the proof of Main Lemma 3.5. Of course, this outline will lack the
majority of the nitty-gritty details, technicalities, that constitute
an integral part of the proof.

\medskip

We consider the one-sided, tubular neighborhoods $U_\delta$ of $J$ with
radius (thickness) $\delta>0$. Throughout the whole proof of the main
lemma the asymptotics of the measures $\mu(X_\delta)$ of certain
(dynamically defined) sets $X_\delta\subset U_\delta$ are studied, as
$\delta\to 0$. We fix a large constant $c_3\gg 1$, and for typical
points $y\in U_\delta\setminus U_{\delta/2}$ (having non-singular
forward orbits and returning to the layer $U_\delta\setminus
U_{\delta/2}$ infinitely many times in the future) we define the
arc-length parametrized curves $\rho_{y,t}(s)$ ($0\le s\le h(y,t)$) in
the following way: $\rho_{y,t}$ emanates from $y$ and it is the curve
inside the manifold $\Sigma_0^t(y)$ with the steepest descent towards
the separating manifold $J$. Here $\Sigma_0^t(y)$ is the inverse image
$S^{-t}\left(\Sigma_t^t(y)\right)$ of the flat, local orthogonal
manifold (flat wave front, see \S2) passing through $y_t=S^t(y)$. The
terminal point $\Pi(y)=\rho_{y,t}\left(h(y,t)\right)$ of the smooth
curve $\rho_{y,t}$ is either

\medskip

(a) on the separating manifold $J$, or

\medskip

(b) on a singularity of order $k_1=k_1(y)$.

\medskip

\noindent
The case (b) is further split in two sub-cases, as follows:

\medskip

(b/1) $k_1(y)<c_3$;

\medskip

(b/2) $c_3\le k_1(y)<\infty$.

\medskip

The set of (typical) points $y\in U_\delta\setminus U_{\delta/2}$ with
property (a) (this is the set $\overline{U}_\delta(\infty)$ in \S3) is
handled by lemmas 3.28 and 3.29, where it is shown that, actually,
$\overline{U}_\delta(\infty)=\emptyset$. Roughly speaking, the reason
of this is the following: For a point
$y\in \overline{U}_\delta(\infty)$ the powers $S^t$ of the flow
exhibit arbitrarily large contractions on the curves $\rho_{y,t}$
(Appendix II), thus the infinitely many returns of $S^t(y)$ to the
layer $U_\delta\setminus U_{\delta/2}$ would "pull up" the other
endpoints $S^t\left(\Pi(y)\right)$ to the region $U_\delta\setminus
J$, consisting entirely of sufficient points, and showing that the
point $\Pi(y)\in J$ itself is sufficient.

The set $\overline{U}_\delta\setminus\overline{U}_\delta(c_3)$ of all
phase points $y\in U_\delta\setminus U_{\delta/2}$ with the property
$k_1(y)<c_3$ are dealt with by Lemma 3.27, where it is shown that

$$
\mu\left(\overline{U}_\delta\setminus\overline{U}_\delta(c_3)\right)=o(\delta),
$$
as $\delta\to 0$. The reason, in rough terms, is that such phase
points must lie at the distance $\le\delta$ from the compact
singularity set

$$
\bigcup\Sb 0\le t\le 2c_3\endSb S^{-t}\left(\Cal S\Cal R^-\right),
$$
and this compact singularity set is transversal to $J$, thus ensuring
the measure estimate
$\mu\left(\overline{U}_\delta\setminus\overline{U}_\delta(c_3)\right)=o(\delta)$.

Finally, the set $F_\delta(c_3)$ of (typical) phase points 
$y\in U_\delta\setminus U_{\delta/2}$ with $c_3\le k_1(y)<\infty$ is
dealt with by lemmas 3.36, 3.37, and Corollary 3.38, where it is shown
that $\mu\left(F_\delta(c_3)\right)\le C\cdot \delta$, with constants
$C$ that can be chosen arbitrarily small by selecting the constant
$c_3\gg 1$ big enough. The ultimate reason of this measure estimate is
the following fact: For every point $y\in F_\delta(c_3)$ the
projection

$$
\tilde{\Pi}(y)=S^{t_{\overline{k}_1(y)}}\in\partial\bold M
$$
(where $t_{\overline{k}_1(y)}$ is the time of the
$\overline{k}_1(y)$-th collision on the forward orbit of $y$) will
have a tubular distance $z_{tub}\left(\tilde{\Pi}(y)\right)\le
C_1\delta$ from the singularity set $\Cal S\Cal R^-\cup \Cal S\Cal R^+$,
where the constant $C_1$ can be made arbitrarily small by choosing the
contraction coefficients of the powers $S^{t_{\overline{k}_1(y)}}$ on
the curves $\rho_{y,t_{\overline{k}_1(y)}}$ arbitrarily small with the
help of the result in Appendix II. The upper masure estimate (inside
the set $\partial\bold M$) of the set of such points 
$\tilde{\Pi}(y)\in\partial\bold M$ (Lemma 2 in [S-Ch(1987)]) finally
yields the required upper bound 
$\mu\left(F_\delta(c_3)\right)\le C\cdot \delta$ with arbitrarily
small positive constants $C$ (if $c_3\gg 1$ is big enough).

The listed measure estimates and the obvious fact

$$
\mu\left(U_\delta\setminus U_{\delta/2}\right)\approx C_2\cdot\delta
$$
(with some constant $C_2>0$, depending only on $J$) show that there must exist a point 
$y\in U_\delta\setminus U_{\delta/2}$ with the property (a) above,
thus ensuring the sufficiency of the point $\Pi(y)\in J$.

\medskip

Finally, in the closing section we complete the inductive proof of ergodicity
(with respect to the number of balls $N$) by utilizing Main Lemma 3.5 and
earlier results from the literature.  Actually, a consequence of the
Main Lemma will be that exceptional $J$-manifolds do not exist, and
this will imply the fact that no distinct, open ergodic components can
coexist.

Appendix I at the end of this paper serves the purpose of making the
reading of the proof of \S3 easier, by providing a chart for the hierarchy of
the selection of several constants playing a role in the proof of 
Main Lemma 3.5.

Appendix II contains a useful (also, potentially useful in the future) uniform
contraction estimate which is exploited in Section 3. Many ideas
of Appendix II originate from N. I. Chernov.

\bigskip \bigskip

\heading
\S2. Prerequisites
\endheading

\bigskip

Consider the $\nu$-dimensional ($\nu\ge2$), standard, flat torus
$\Bbb T^\nu=\Bbb R^\nu/\Bbb Z^\nu$ as the vessel containing 
$N$ ($\ge2$) hard balls (spheres) $B_1,\dots,B_N$ with positive masses 
$m_1,\dots,m_N$ and (just for simplicity) common radius $r>0$. We always
assume that the radius $r>0$ is not too big, so
that even the interior of the arising configuration space $\bold Q$ (or, 
equivalently, the phase space) is connected. Denote the center of the ball
$B_i$ by $q_i\in\Bbb T^\nu$, and let $v_i=\dot q_i$ be the velocity of the
$i$-th particle. We investigate the uniform motion of the balls
$B_1,\dots,B_N$ inside the container $\Bbb T^\nu$ with half a unit of total 
kinetic energy: $E=\dfrac{1}{2}\sum_{i=1}^N m_i||v_i||^2=\dfrac{1}{2}$.
We assume that the collisions between balls are perfectly elastic. Since
--- beside the kinetic energy $E$ --- the total momentum
$I=\sum_{i=1}^N m_iv_i\in\Bbb R^\nu$ is also a trivial first integral of the
motion, we make the standard reduction $I=0$. Due to the apparent translation
invariance of the arising dynamical system, we factorize the configuration
space with respect to uniform spatial translations as follows:
$(q_1,\dots,q_N)\sim(q_1+a,\dots,q_N+a)$ for all translation vectors
$a\in\Bbb T^\nu$. The configuration space $\bold Q$ of the arising flow
is then the factor torus
$\left(\left(\Bbb T^\nu\right)^N/\sim\right)\cong\Bbb T^{\nu(N-1)}$
minus the cylinders

$$
C_{i,j}=\left\{(q_1,\dots,q_N)\in\Bbb T^{\nu(N-1)}\colon\;
\text{dist}(q_i,q_j)<2r \right\}
$$
($1\le i<j\le N$) corresponding to the forbidden overlap between the $i$-th
and $j$-th spheres. Then it is easy to see that the compound 
configuration point

$$
q=(q_1,\dots,q_N)\in\bold Q=\Bbb T^{\nu(N-1)}\setminus
\bigcup_{1\le i<j\le N}C_{i,j}
$$
moves in $\bold Q$ uniformly with unit speed and bounces back from the
boundaries $\partial C_{i,j}$ of the cylinders $C_{i,j}$ according to the
classical law of geometric optics: the angle of reflection equals the angle of
incidence. More precisely: the post-collision velocity $v^+$ can be obtained
from the pre-collision velocity $v^-$ by the orthogonal reflection across the
tangent hyperplane of the boundary $\partial\bold Q$ at the point of collision.
Here we must emphasize that the phrase ``orthogonal'' should be understood 
with respect to the natural Riemannian metric (the kinetic energy)
$||dq||^2=\sum_{i=1}^N m_i||dq_i||^2$ in the configuration space $\bold Q$.
For the normalized Liouville measure $\mu$ of the arising flow
$\{S^t\}$ we obviously have $d\mu=\text{const}\cdot dq\cdot dv$, where
$dq$ is the Riemannian volume in $\bold Q$ induced by the above metric, 
and $dv$ is the surface measure (determined by the restriction of the
Riemannian metric above) on the unit sphere of compound velocities

$$
\Bbb S^{\nu(N-1)-1}=\left\{(v_1,\dots,v_N)\in\left(\Bbb R^\nu\right)^N\colon\;
\sum_{i=1}^N m_iv_i=0 \text{ and } \sum_{i=1}^N m_i||v_i||^2=1 \right\}.
$$
The phase space $\bold M$ of the flow $\{S^t\}$ is the unit tangent bundle
$\bold Q\times\Bbb S^{d-1}$ of the configuration space $\bold Q$. (We will 
always use the shorthand notation $d=\nu(N-1)$ for the dimension of the 
billiard table $\bold Q$.) We must, however, note here that at the boundary
$\partial\bold Q$ of $\bold Q$ one has to glue together the pre-collision and
post-collision velocities in order to form the phase space $\bold M$, so
$\bold M$ is equal to the unit tangent bundle $\bold Q\times\Bbb S^{d-1}$
modulo this identification.

A bit more detailed definition of hard ball systems with arbitrary masses,
as well as their role in the family of cylindric billiards, can be found in
\S4 of [S-Sz(2000)] and in \S1 of [S-Sz(1999)]. We denote the
arising flow by $\flow$.

In the late 1970s Sinai [Sin(1979)] developed a powerful, three-step strategy
for proving the (hyperbolic) ergodicity of hard ball systems. This strategy
was later implemented in a series of papers [K-S-Sz(1989)], [K-S-Sz(1990)-I],
[K-S-Sz(1991)], and [K-S-Sz(1992)]. First of all, these proofs are inductions
on the number $N$ of balls involved in the problem. Secondly, the induction
step itself consists of the following three major steps:

\medskip

\subheading{Step I} To prove that every non-singular (i. e. smooth)
trajectory segment $\traj$ with a ``combinatorially rich'' (in a well defined
sense of Definition 3.28 of [Sim(2002)]) symbolic collision sequence is
automatically sufficient (or, in other words, ``geometrically hyperbolic'',
see below in this section), provided that the phase point $x_0$ does not
belong to a countable union $J$ of smooth sub-manifolds with codimension at
least two. (Containing the exceptional phase points.)

The exceptional set $J$ featuring this result is negligible in our dynamical
considerations --- it is a so called slim set. For the basic properties of
slim sets, see again below in this section.

\medskip

\subheading{Step II} Assume the induction hypothesis, i. e. that all hard
ball systems with $N'$ balls ($2\le N'<N$) are (hyperbolic and) ergodic.
Prove that there exists a slim set $E\subset\bold M$ with the following
property: For every phase point $x_0\in\bold M\setminus E$ the entire
trajectory $S^{\Bbb R}x_0$ contains at most one singularity and its symbolic
collision sequence is combinatorially rich in the sense of Definition 3.28 of
[Sim(2002)], just as required by the result of Step I.

\medskip

\subheading{Step III} By using again the induction hypothesis, prove that
almost every singular trajectory is sufficient in the time interval
$(t_0,+\infty)$, where $t_0$ is the time moment of the singular reflection.
(Here the phrase ``almost every'' refers to the volume defined by the induced
Riemannian metric on the singularity manifolds.)

We note here that the almost sure sufficiency of the singular trajectories
(featuring Step III) is an essential condition for the proof of the celebrated
Theorem on Local Ergodicity for semi-dispersive billiards proved by Chernov 
and Sinai [S-Ch(1987)]. Under this assumption, the result of 
Chernov and Sinai states that in any semi-dispersive billiard system a 
suitable, open neighborhood $U_0$ of any sufficient phase point 
$x_0\in\bold M$ (with at most one singularity on its trajectory) belongs
to a single ergodic component of the billiard flow $\flow$. 

A few years ago B\'alint, Chernov, Sz\'asz, and T\'oth [B-Ch-Sz-T(2002)]
discovered that, in addition, the algebraic nature of the scatterers needs
to be assumed, in order for the proof of this result to work. Fortunately,
systems of hard balls are, by nature, automatically algebraic.

In an inductive proof of ergodicity, steps I and II together ensure that
there exists an arc-wise connected set
$C\subset\bold M$ with full measure, such that every phase point $x_0\in C$
is sufficient with at most one singularity on its trajectory. Then the cited
Theorem on Local Ergodicity (now taking advantage of the result of Step III)
states that for every phase point $x_0\in C$ an open neighborhood $U_0$ of
$x_0$ belongs to one ergodic component of the flow. Finally, the connectedness
of the set $C$ and $\mu(\bold M\setminus C)=0$ imply that the flow
$\flow$ (now with $N$ balls) is indeed ergodic, and actually fully hyperbolic,
as well.

\medskip

The generator subspace $A_{i,j}\subset \Bbb R^{\nu N}$ ($1\le i<j\le N$)
of the cylinder $C_{i,j}$ (describing the collisions between the $i$-th and
$j$-th balls) is given by the equation
$$
A_{i,j}=\left\{(q_1,\dots,q_N)\in\left(\Bbb R^\nu\right)^N\colon\;
q_i=q_j \right\},
\tag 2.1
$$
see (4.3) in [S-Sz(2000)]. Its ortho-complement 
$L_{i,j}\subset\Bbb R^{\nu N}$ is then defined by the equation
$$
L_{i,j}=\left\{(q_1,\dots,q_N)\in\left(\Bbb R^\nu\right)^N\colon\;
q_k=0 \text{ for } k\ne i,j, \text{ and } m_iq_i+m_jq_j=0 \right\},
\tag 2.2
$$
see (4.4) in [S-Sz(2000)].
Easy calculation shows that the cylinder $C_{i,j}$ 
(describing the overlap of the $i$-th and $j$-th balls)
is indeed spherical and the radius of its base sphere is equal to
$r_{i,j}=2r\sqrt{\frac{m_im_j}{m_i+m_j}}$, see \S 4, especially formula
(4.6) in [S-Sz(2000)].

The structure lattice $\Cal L\subset\Bbb R^{\nu N}$ is clearly the
lattice $\Cal L=\left(\Bbb Z^{\nu}\right)^N=\Bbb Z^{N\nu}$. 

Due to the presence of an additional invariant quantity
$I=\sum_{i=1}^N m_iv_i$, one usually makes the reduction
$\sum_{i=1}^N m_iv_i=0$ and, correspondingly, factorizes the configuration
space with respect to uniform spatial translations:

$$
(q_1,\dots,q_N)\sim(q_1+a,\dots,q_N+a), \quad a\in\Bbb T^\nu.
\tag 2.3
$$
The natural, common tangent space of this reduced configuration space is

$$
\Cal Z=\left\{(v_1,\dots,v_N)\in\left(\Bbb R^\nu\right)^N\colon\;
\sum_{i=1}^N m_iv_i=0\right\}=\left(\bigcap_{i<j}A_{i,j}
\right)^\perp=\left(\Cal A\right)^\perp
\tag 2.4
$$
supplied with the inner product 

$$
\langle v,\,v'\rangle=\sum_{i=1}^N m_i\langle v_i,\,v_i'\rangle,
$$
see also (4.1) and (4.2) in [S-Sz(2000)]. 

\bigskip

\subheading{Collision graphs} Let $S^{[a,b]}x$ be a nonsingular, finite
trajectory segment with the collisions $\sigma_1,\dots,\sigma_n$
listed in time order. 
(Each $\sigma_k$ is an unordered pair $(i,j)$ of different labels
$i,j\in\{1,2,\dots,N\}$.) The graph $\Cal G=(\Cal V,\Cal E)$ with vertex set
$\Cal V=\{1,2,\dots,N\}$ and set of edges $\Cal E=\{\sigma_1,\dots,\sigma_n\}$ 
is called the {\it collision graph}
of the orbit segment $S^{[a,b]}x$. For a given positive number $C$, the
collision graph $\Cal G=(\Cal V,\Cal E)$ of the orbit segment $S^{[a,b]}x$
will be called {\it $C$-rich} if $\Cal G$ contains at least $C$ connected,
consecutive (i. e. following one after the other in time, according to the
time-ordering given by the trajectory segment $S^{[a,b]}x$) subgraphs. 

\bigskip

\subheading{Singularities and Trajectory Branches} There are two types of
singularities of the billiard flow: tangential (or gliding) and double
collisions of balls. The first means that two balls collide in such a way that
their relative velocity vector is parallel to their common tangent hyperplane
at the point of contact. In this case no momentum is exchanged, and on the two
sides of this singularity the flow behaves differently: on one side the two
balls collide with each other, on the other side they fly by each other
without interaction.

The second type of singularity, a double collision means that two pairs of
particles, $(i,j)$ and $(j,k)$ ($i$, $j$, and $k$ are three different labels)
are to collide exactly at the same time. (The case when the two interacting
pairs do not have a common particle $j$ may be disregarded, since this does
not give rise to non-differentiability of the flow.)

We are going to briefly describe the discontinuity of the flow
$\{S^t\}$ caused by a double collisions at time $t_0$.
Assume first that the pre--collision velocities of the particles are given.
What can we say about the possible post--collision velocities? Let us perturb
the pre--collision phase point (at time $t_0-0$) infinitesimally, so that the
collisions at $\sim t_0$ occur at infinitesimally different moments. By
applying the collision laws to the arising finite sequence of collisions, we
see that the post-collision velocities are fully determined by the 
time-ordered list of the arising collisions. Therefore, the collection of all
possible time-ordered lists of these collisions gives rise to a finite family 
of continuations of the trajectory beyond $t_0$. They are called the
trajectory branches. It is quite clear that similar statements can be
said regarding the evolution of a trajectory through a double collision
in reverse time. Furthermore, it is also obvious that for any given
phase point $x_0\in\bold M$ there are two, $\omega$-high trees
$\Cal T_+$ and $\Cal T_-$ such that $\Cal T_+$ ($\Cal T_-$) describes all the
possible continuations of the positive (negative) trajectory
$S^{[0,\infty)}x_0$ ($S^{(-\infty,0]}x_0$). (For the definitions of trees and
for some of their applications to billiards, cf. the beginning of \S 5
in [K-S-Sz(1992)].) It is also clear that all possible continuations
(branches) of the whole trajectory $S^{(-\infty,\infty)}x_0$ can be uniquely
described by all pairs $(B_-,B_+)$ of infinite branches of
the trees $\Cal T_-$ and $\Cal T_+$ ($B_-\subset\Cal T_-, B_+\subset
\Cal T_+$).

Finally, we note that the trajectory of the phase point $x_0$ has exactly two
branches, provided that $S^tx_0$ hits a singularity for a single value
$t=t_0$, and the phase point $S^{t_0}x_0$ does not lie on the intersection
of more than one singularity manifolds. In this case we say that the 
trajectory of $x_0$ has a ``simple singularity''.

Other singularities (phase points lying on the intersections of two or more
singularity manifolds of codimension $1$) can be disregarded in our studies,
since these points lie on a countable collection of codimension-two, smooth
submanifolds of the phase space, and such a special type of slim set can
indeed be safely discarded in our proof, see the part {\it Slim sets} later in
this section.

\bigskip

\subheading{Neutral Subspaces, Advance, and Sufficiency}
Consider a nonsingular trajectory segment $S^{[a,b]}x$.
Suppose that $a$ and $b$ are not moments of collision.

\medskip

\proclaim{Definition 2.5} The neutral space $\Cal N_0(S^{[a,b]}x)$
of the trajectory segment $S^{[a,b]}x$ at time zero ($a<0<b$)\ is
defined by the following formula:

$$
\aligned
&\Cal N_0(S^{[a,b]}x)=\big \{W\in\Cal Z\colon\;\exists (\delta>0) \;
\text{ such that } \; \forall \alpha \in (-\delta,\delta) \\
&V\left(S^a\left(Q(x)+\alpha W,V(x)\right)\right)=V(S^ax)\text{ and }
V\left(S^b\left(Q(x)+\alpha W,V(x)\right)\right)=V(S^bx)\big\}.
\endaligned
$$
\endproclaim

\noindent
($\Cal Z$ is the common tangent space $\Cal T_q\bold Q$ of the parallelizable
manifold $\bold Q$ at any of its points $q$, while $V(x)$ is the velocity
component of the phase point $x=\left(Q(x),\,V(x)\right)$.)

It is known (see (3) in \S 3 of [S-Ch (1987)]) that
$\Cal N_0(S^{[a,b]}x)$ is a linear subspace of $\Cal Z$ indeed, and
$V(x)\in \Cal N_0(S^{[a,b]}x)$. The neutral space $\Cal N_t(S^{[a,b]}x)$
of the segment $S^{[a,b]}x$ at time $t\in [a,b]$ is defined as follows:

$$
\Cal N_t(S^{[a,b]}x)=\Cal N_0\left(S^{[a-t,b-t]}(S^tx)\right).
$$
It is clear that the neutral space $\Cal N_t(S^{[a,b]}x)$ can be
canonically
identified with $\Cal N_0(S^{[a,b]}x)$ by the usual identification of the
tangent spaces of $\bold Q$ along the trajectory $S^{(-\infty,\infty)}x$
(see, for instance, \S 2 of [K-S-Sz(1990)-I]).

Our next  definition is  that of  the advance. Consider a
non-singular orbit segment $S^{[a,b]}x$ with the symbolic collision sequence
$\Sigma=(\sigma_1,\dots,\sigma_n)$, meaning that $S^{[a,b]}x$
has exactly $n$ collisions with $\partial\bold Q$, and the $i$-th collision
($1\le i\le n$) takes place at the boundary of the cylinder $C_{\sigma_i}$.
For $x=(Q,V)\in\bold M$ and $W\in\Cal Z$, $\Vert W\Vert$ sufficiently small, 
denote $T_W(Q,V):=(Q+W,V)$.

\proclaim{Definition 2.6}
For any $1\le k\le n$ and $t\in[a,b]$, the advance

$$
\alpha_k=\alpha(\sigma_k)\colon\;\Cal N_t(S^{[a,b]}x)\rightarrow\Bbb R
$$
of the collision $\sigma_k$ is the unique linear extension of the linear 
functional $\alpha_k=\alpha(\sigma_k)$
defined in a sufficiently small neighborhood of the origin of 
$\Cal N_t(S^{[a,b]}x)$ in the following way:
$$
\alpha(\sigma_k)(W):= t_k(x)-t_k(S^{-t}T_WS^tx).
$$
\endproclaim

Here $t_k=t_k(x)$ is the time of the $k$-th collision $\sigma_k$ on
the trajectory of $x$ after time $t=a$. The above formula and the notion of
the advance functional 

$$
\alpha_k=\alpha(\sigma_k):\; \Cal N_t\left(S^{[a,b]}x\right)\to\Bbb R
$$
has two important features:

\medskip

(i) If the spatial translation $(Q,V)\mapsto(Q+W,V)$ 
($W\in\Cal N_t\left(S^{[a,b]}x\right)$) is carried out at time
$t$, then $t_k$ changes linearly in $W$, and it takes place just 
$\alpha_k(W)$ units of time earlier. (This is why it is called ``advance''.)

\medskip

(ii) If the considered reference time $t$ is somewhere between $t_{k-1}$
and $t_k$, then the neutrality of $W$ with respect to $\sigma_k$ 
precisely means that

$$
W-\alpha_k(W)\cdot V(x)\in A_{\sigma_k},
$$
i. e. a neutral (with respect to the collision $\sigma_k$) spatial 
translation $W$ with the advance $\alpha_k(W)=0$ means that the vector $W$ 
belongs to the generator space
$A_{\sigma_k}$ of the cylinder $C_{\sigma_k}$.

It is now time to bring up the basic notion of sufficiency 
(or, sometimes it is also called geometric hyperbolicity) of a
trajectory (segment). This is the utmost important necessary condition for
the proof of the Theorem on Local Ergodicity for semi-dispersive billiards,
[S-Ch(1987)].

\medskip

\proclaim{Definition 2.7}
\roster
\item
The nonsingular trajectory segment $S^{[a,b]}x$ ($a$ and $b$ are supposed not
to be moments of collision) is said to be sufficient if and only if
the dimension of $\Cal N_t(S^{[a,b]}x)$ ($t\in [a,b]$) is minimal, i.e.
$\text{dim}\ \Cal N_t(S^{[a,b]}x)=1$.
\item
The trajectory segment $S^{[a,b]}x$ containing exactly one singularity (a so 
called ``simple singularity'', see above) is said to be sufficient if 
and only if both branches of this trajectory segment are sufficient.
\endroster
\endproclaim

\medskip

\proclaim{Definition 2.8}
The phase point $x\in\bold M$ with at most one (simple) singularity is said
to be sufficient if and only if its whole trajectory $S^{(-\infty,\infty)}x$
is sufficient, which means, by definition, that some of its bounded
segments $S^{[a,b]}x$ are sufficient.
\endproclaim

\medskip

\subheading{Note} In this paper the phrase "trajectory (segment) with at most
one singularity" always means that the sole singularity of the trajectory
(segment), if exists, is simple.

\medskip

In the case of an orbit $S^{(-\infty,\infty)}x$ with at most one
singularity, sufficiency means that both branches of
$S^{(-\infty,\infty)}x$ are sufficient.

\bigskip

\subheading{No accumulation (of collisions) in finite time} 
By the results of Vaserstein [V(1979)], Galperin [G(1981)] and
Burago-Ferleger-Kononenko [B-F-K(1998)], in any semi-dis\-per\-sive 
billiard flow there can only be finitely many 
collisions in finite time intervals, see Theorem 1 in [B-F-K(1998)]. 
Thus, the dynamics is well defined as long as the trajectory does not hit 
more than one boundary components at the same time.

\bigskip

\subheading{Slim sets} We are going to summarize the basic properties of
codimension-two subsets $A$ of a connected, smooth manifold $M$ with a
possible boundary and corners.  Since these subsets $A$ are just those
negligible in our dynamical discussions, we shall call them slim. As to a
broader exposition of the issues, see [E(1978)] or \S2 of [K-S-Sz(1991)].

Note that the dimension $\dim A$ of a separable metric space $A$ is one of the
three classical notions of topological dimension: the covering 
(\v Cech-Lebesgue), the small inductive (Menger-Urysohn), or the large 
inductive (Brouwer-\v Cech) dimension. As it is known from general
topology, all of them are the same for separable metric spaces, 
see [E(1978)].

\medskip

\subheading{Definition 2.9}
A subset $A$ of $M$ is called slim if and only if $A$ can be covered by a 
countable family of codimension-two (i. e. at least two) closed sets of
$\mu$--measure zero, where $\mu$ is any smooth measure on $M$. (Cf.
Definition 2.12 of [K-S-Sz(1991)].)

\medskip

\subheading{Property 2.10} The  collection of all slim subsets of $M$ is a
$\sigma$-ideal, that is, countable unions of slim sets and arbitrary
subsets of slim sets are also slim.

\medskip

\subheading{Proposition 2.11. (Locality)}
A subset $A\subset M$ is slim if and only if for
every $x\in A$ there exists an open neighborhood $U$ of $x$ in $M$ such that 
$U\cap A$ is slim. (Cf. Lemma 2.14 of [K-S-Sz(1991)].)

\medskip

\subheading{Property 2.12} A closed subset $A\subset M$ is slim if and only
if $\mu(A)=0$ and $\dim A\le\dim M-2$.

\medskip

\subheading{Property 2.13. (Integrability)}
If $A\subset M_1\times M_2$ is a closed subset of the product of two smooth,
connected manifolds with possible boundaries and corners, and for every 
$x\in M_1$ the set

$$
A_x=\{ y\in M_2\colon\; (x,y)\in A\}
$$
is slim in $M_2$, then $A$ is slim in $M_1\times M_2$.

\medskip

The following propositions characterize the codimension-one and 
codimension-two sets.

\medskip

\subheading{Proposition 2.14}
For any closed subset $S\subset M$ the following three conditions are 
equivalent:

\roster

\item"{(i)}" $\dim S\le\dim M-2$;

\item"{(ii)}"  $\text{int}S=\emptyset$ and for every open connected set 
$G\subset M$ the difference set $G\setminus S$ is also connected;

\item"{(iii)}" $\text{int}S=\emptyset$ and for every point $x\in M$ and for any
open neighborhood $V$ of $x$ in $M$ there exists a smaller open neighborhood
$W\subset V$ of the point $x$ such that for every pair of points 
$y,z\in W\setminus S$ there is a continuous curve $\gamma$ in the set 
$V\setminus S$ connecting the points $y$ and $z$.

\endroster

\noindent
(See Theorem 1.8.13 and Problem 1.8.E of [E(1978)].)

\medskip

\subheading{Proposition 2.15} For any subset $S\subset M$ the condition 
$\dim S\le\dim M-1$ is equivalent to $\text{int}S=\emptyset$.
(See Theorem 1.8.10 of [E(1978)].)

\medskip

We recall an elementary, but important lemma (Lemma 4.15 of [K-S-Sz(1991)]).
Let $\Delta_2$ be the set of phase points 
$x\in\bold M\setminus\partial\bold M$ such that the trajectory 
$S^{(-\infty,\infty)}x$ has more than one singularities (or, its only 
singularity is not simple).

\subheading{Proposition 2.16}
The set $\Delta_2$ is a countable union of codimension-two
smooth submanifolds of $M$ and, being such, is slim.

\medskip

The next lemma establishes the most important property of slim sets which
gives us the fundamental geometric tool to connect the open ergodic components
of billiard flows.

\medskip

\subheading{Proposition 2.17}
If $M$ is connected, then the complement $M\setminus A$ of a slim $F_\sigma$ 
set $A\subset M$ is an arc-wise connected ($G_\delta$) set of full measure. 
(See Property 3 of \S 4.1 in [K-S-Sz(1989)]. The  $F_\sigma$ sets are, 
by definition, the countable unions of closed sets, while the $G_\delta$ sets
are the countable intersections of open sets.)

\medskip

\subheading{The subsets $\bold M^0$ and $\bold M^\#$} Denote by
$\bold M^\#$ the set of all phase points $x\in\bold M$ for which the
trajectory of $x$ encounters infinitely many non-tangential collisions
in both time directions. The trajectories of the points 
$x\in\bold M\setminus\bold M^\#$ are lines: the motion is linear and uniform,
see the appendix of [Sz(1994)]. It is proven in lemmas A.2.1 and A.2.2
of [Sz(1994)] that the closed set $\bold M\setminus\bold M^\#$ is a finite
union of hyperplanes. It is also proven in [Sz(1994)] that, locally, the two
sides of a hyper-planar component of $\bold M\setminus\bold M^\#$ can be
connected by a positively measured beam of trajectories, hence, from the point
of view of ergodicity, in this paper it is enough to show that the connected
components of $\bold M^\#$ entirely belong to one ergodic component. This is
what we are going to do in this paper.

Denote by $\bold M^0$ the set of all phase points $x\in\bold M^\#$ the 
trajectory of which does not hit any singularity, and use the notation
$\bold M^1$ for the set of all phase points $x\in\bold M^\#$ whose orbit
contains exactly one, simple singularity. According to Proposition 2.16,
the set $\bold M^\#\setminus(\bold M^0\cup\bold M^1)$ is a countable union of
smooth, codimension-two ($\ge2$) submanifolds of $\bold M$, and, therefore,
this set may be discarded in our study of ergodicity, please see also the
properties of slim sets above. Thus, we will restrict our attention to the
phase points $x\in\bold M^0\cup\bold M^1$.

\medskip

\subheading{The ``Chernov-Sinai Ansatz''} An essential precondition for
the Theorem on Local Ergodicity by Chernov and Sinai [S-Ch(1987)] is the
so called ``Chernov-Sinai Ansatz'' which we are going to formulate below.
Denote by $\Cal S\Cal R^+\subset\partial\bold M$ the set of all phase points
$x_0=(q_0,v_0)\in\partial\bold M$ corresponding to singular reflections
(a tangential or a double collision at time zero) supplied with the 
post-collision (outgoing) velocity $v_0$. It is well known that
$\Cal S\Cal R^+$ is a compact cell complex with dimension
$2d-3=\text{dim}\bold M-2$. It is also known (see Lemma 4.1 in 
[K-S-Sz(1990)-I], in conjunction with Proposition 2.16 above)
that for $\nu_1$-almost every phase point $x_0\in\Cal S\Cal R^+$ the forward 
orbit $S^{(0,\infty)}x_0$ does not hit any further singularity. 
(Here $\nu_1$ is the Riemannian volume of $\Cal S\Cal R^+$ induced by the 
restriction of the natural Riemannian metric of $\bold M$.)
The Chernov-Sinai Ansatz postulates that for $\nu_1$-almost every 
$x_0\in\Cal S\Cal R^+$ the forward orbit $S^{(0,\infty)}x_0$ is sufficient 
(geometrically hyperbolic).

\medskip

\subheading{The Theorem on Local Ergodicity} The Theorem on Local 
Ergodicity for semi-dispersive billiards 
(Theorem 5 of [S-Ch(1987)]) claims the following: 
Let $\flow$ be a semi-dispersive billiard flow with the property
that the smooth components of the boundary $\partial\bold Q$ of the 
configuration space are algebraic hyper-surfaces. (The cylindric billiards
automatically fulfill this algebraicity condition.) Assume -- further --
that the Chernov-Sinai Ansatz holds true, and a phase point 
$x_0\in\left(\bold M^0\cup\bold M^1\right)\setminus\partial\bold M$ 
is sufficient.

\medskip

Then some open neighborhood $U_0\subset\bold M$ of $x_0$ belongs to a single
ergodic component of the flow $\flow$. (Modulo the zero sets, of course.)

\bigskip \bigskip

\heading
\S3. The Exceptional $J$-Manifolds \\
(The asymptotic measure estimates)
\endheading

\bigskip \bigskip

First of all, we define the fundamental object for the proof of our theorem.

\medskip

\subheading{Definition 3.1} A smooth submanifold $J\subset\text{int}\bold M$
of the interior of the phase space $\bold M$ is called an {\it exceptional
$J$-manifold} (or simply an exceptional manifold) with a negative Lyapunov
function $Q$ if

\medskip

(1) $\text{dim}J=2d-2$ ($=\text{dim}\bold M-1$);

\medskip

(2) the pair of manifolds $(\overline{J},\,\partial J)$ is diffeomorphic 
to the standard pair 

$$
\left(B^{2d-2},\,\Bbb S^{2d-3}\right)=\left(B^{2d-2},\,
\partial B^{2d-2}\right),
$$
where $B^{2d-2}$ is the closed unit ball of $\Bbb R^{2d-2}$;

\medskip

(3) $J$ is locally flow-invariant, i. e. $\forall x\in J$ 
$\exists\,a(x),\,b(x)$, $a(x)<0<b(x)$, such that $S^tx\in J$ for all $t$ with
$a(x)<t<b(x)$, and $S^{a(x)}x\in\partial J$, $S^{b(x)}x\in\partial J$;

\medskip

(4) the manifold $J$ has some thin, open, tubular neighborhood $\tilde U_0$ in
$\text{int}\bold M$, and there exists a number $T>0$ such that

(i) $S^T\left(\tilde{U}_0\right)\cap\partial\bold M=\emptyset$, and
all orbit segments $S^{[0,T]}x$ ($x\in\tilde U_0$) are non-singular, hence
they share the same symbolic collision sequence $\Sigma$;

(ii) $\forall x\in\tilde U_0$ the orbit segment $S^{[0,T]}x$ is sufficient
if and only if $x\not\in J$;

\medskip

(5) $\forall x\in J$ we have $Q(n(x)):=\langle z(x),\,w(x)\rangle\le -c_1<0$
for a unit normal vector field $n(x)=(z(x),\,w(x))$ of $J$ with a fixed
constant $c_1>0$;

\medskip

(6) the set $W$ of phase points $x\in J$ never again returning to $J$ (After 
first leaving it, of course. Keep in mind that $J$ is locally flow-invariant)
has relative measure greater than $1-10^{-8}$ in $J$, i. e.
$\dfrac{\mu_1(W)}{\mu_1(J)}>1-10^{-8}$, where $\mu_1$ is the hypersurface 
measure of the smooth manifold $J$.

\medskip

\subheading{Remark} The above definition is, by nature, fairly technical, thus
a short commenting of it is due here. Once we make the induction hypothesis,
i. e.  we assume that the (hyperbolic) ergodicity and the Chernov-Sinai Ansatz
hold true for any hard ball systems with less than $N$ balls (regardless of
the masses $m_i$ and the radius $r$), the only way for two distinct ergodic
components to co-exist is when they are separated by an exceptional manifold
$J$ described in the above definition. This is proved in \S4 below.

\medskip

We begin with an important proposition on the structure of forward orbits
$S^{[0,\infty)}x$ for $x\in J$.

\medskip

\subheading{Proposition 3.2} For $\mu_1$-almost every $x\in J$
the forward orbit $S^{[0,\infty)}x$ is non-singular.

\medskip

\subheading{Proof} According to Proposition 7.12 of [Sim(2003)], the set

$$
J\cap\left[\bigcup\Sb t>0\endSb S^{-t}\left(\Cal S\Cal R^-\right)\right]
$$
of forward singular points $x\in J$ is a countable union of smooth, proper
submanifolds of $J$, hence it has $\mu_1$-measure zero. \qed

\medskip

In the future we will need

\medskip

\subheading{Lemma 3.3} The concave, local orthogonal manifolds $\Sigma(y)$
passing through points $y\in J$ are uniformly transversal to $J$.

\medskip

\subheading{Note} A local orthogonal manifold $\Sigma\subset\text{int}\bold M$
is obtained from a codimension-one, smooth submanifold $\Sigma_1$ of 
$\text{int}\bold Q$ by supplying $\Sigma_1$ with a selected field of unit
normal vectors as velocities. $\Sigma$ is said to be concave if the second
fundamental form of $\Sigma_1$ (with respect to the selected field of
normal vectors)
is negative semi-definite at every point of $\Sigma_1$. Similarly, the 
convexity of $\Sigma$ requires positive semi-definiteness here, 
see also \S2 of [K-S-Sz(1990)-I].

\medskip

\subheading{Proof} We will only prove the transversality. It will
be clear from the uniformity of the estimates used in the proof that the
claimed transversalities are actually uniform across $J$.

Assume, to the contrary of the transversality, that a concave, local orthogonal
manifold $\Sigma(y)$ is tangent to $J$ at some $y\in J$. Let 
$(\delta q,\,B\delta q)$ be any vector of $\Cal T_y\bold M$ tangent to
$\Sigma(y)$ at $y$. Here $B\le 0$ is the second fundamental form of the
projection $q\left(\Sigma(y)\right)=\Sigma_1(y)$ of $\Sigma(y)$ at the point
$q=q(y)$. The assumed tangency means that
$\langle\delta q,\,z\rangle+\langle B\delta q,\,w\rangle=0$, where
$n(y)=\left(z(y),\,w(y)\right)=(z,w)$ is the unit normal vector of $J$ at $y$.
We get that $\langle\delta q,\,z+Bw\rangle=0$ for any vector 
$\delta q\in v(y)^\perp$. We note that the components $z$ and $w$ of $n$ are
necessarily orthogonal to the velocity $v(y)$, because the manifold $J$ is
locally flow-invariant and the velocity is normalized to $1$ in the phase 
space $\bold M$. The last equation means that $z=-Bw$, thus
$Q(n(y))=\langle z,\,w\rangle=\langle -Bw,\,w\rangle\ge0$, contradicting to
the assumption $Q(n(y))\le -c_1$ of (5) in 3.1. This finishes the proof of
the lemma. \qed

\medskip

In order to formulate the main result of this section, we need to define two
important subsets of $J$.

\medskip

\subheading{Definition 3.4} Let

$$
A=\left\{x\in J\big|\; S^{[0,\infty)}x\text{ is nonsingular and }
\text{dim}\Cal N_0\left(S^{[0,\infty)}x\right)=1 \right\},
$$

$$
B=\left\{x\in J\big|\; S^{[0,\infty)}x\text{ is nonsingular and }
\text{dim}\Cal N_0\left(S^{[0,\infty)}x\right)>1\right\}.
$$

The two Borel subsets $A$ and $B$ of $J$ are disjoint and, according to 
Proposition 3.2 above, their union $A\cup B$ has full $\mu_1$-measure in $J$.

The anticipated main result of this section is

\medskip

\subheading{Main Lemma 3.5} Use all of the above definitions and notations.
We claim that $A\ne\emptyset$.

\medskip

\subheading{Proof} The proof will be a proof by contradiction, and it will be
subdivided into several lemmas. Thus, from now on, we assume that 
$A=\emptyset$.

First, select and fix a non-periodic point (a ``base point'') $x_0\in B$.
Following the fundamental construction of local stable invariant manifolds
[S-Ch(1987)] (see also \S5 of [K-S-Sz(1990)-I]), for any $y\in\bold M$ and any $t>0$
we define the concave, local orthogonal manifolds

$$
\aligned
&\Sigma_t^t(y)=SC_{y_t}\left(\left\{(q,\,v(y_t))\in\bold M\big|\; 
q-q(y_t)\perp v(y_t)\right\}\setminus(\Cal S_1\cup\Cal S_{-1})\right), \\
&\Sigma_0^t(y)=SC_y\left[S^{-t}\Sigma_t^t(y)\right],
\endaligned
\tag 3.6
$$
where $\Cal S_1:=\left\{x\in\bold M\big|\; Tx\in\Cal S\Cal R^-\right\}$
(the set of phase points on singularities of order $1$),
$\Cal S_{-1}:=\left\{x\in\bold M\big|\; -x\in\Cal S_1\right\}$
(the set of phase points on singularities of order $-1$),
$y_t=S^ty$, and $SC_y(\,.\,)$ stands for taking the smooth component of 
the given set that contains the point $y$. The local, stable invariant
manifold $\gamma^{s}(y)$ of $y$ is known to be a superset of the
$C^2$-limiting manifold $\lim_{t\to\infty}\Sigma_0^t(y)$.

\medskip

For any $y\in\bold M$ we use the traditional notations

$$
\aligned
&\tau(y)=\min\left\{t>0\big|\; S^ty\in\partial\bold M\right\}, \\
&T(y)=S^{\tau(y)}y
\endaligned
\tag 3.7
$$
for the first hitting of the collision space $\partial\bold M$. The first
return map (Poincar\'e section, collision map) 
$T:\,\partial\bold M\to\partial\bold M$ (the restriction of the above $T$
to $\partial\bold M$) is known to preserve the finite measure $\nu$ that can
be obtained from the Liouville measure $\mu$ by projecting the latter one onto
$\partial\bold M$ along the flow. Following 4. of [K-S-Sz(1990)-II], for any
point $y\in\text{int}\bold M$ (with $\tau(y)<\infty$, $\tau(-y)<\infty$, where
$-y=(q,-v)$ for $y=(q,v)$) we denote by $z_{tub}(y)$ the supremum of all radii
$\rho>0$ of tubular neighborhoods $V_\rho$ of the projected segment

$$
q\left(\left\{S^ty\big|\; -\tau(-y)\le t \le\tau(y)\right\}\right)
\subset\bold Q
$$
for which even the closure of the set

$$
\left\{(q,\,v(y))\in\bold M\big|\; q\in V_\rho\right\}
$$
does not intersect the set $\Cal S_1\cup\Cal S_{-1}$.
We remind the reader that both Lemma 2 of [S-Ch(1987)] and Lemma 4.10 of
[K-S-Sz(1990)-I] use this tubular distance function $z_{tub}(\,.\,)$ (despite
the notation $z(\,.\,)$ in those papers), see the important note 4. in
[K-S-Sz(1990)-II].

On all the arising local orthogonal manifolds, appearing in the proof, we will
always use the so called $\delta q$-metric to measure distances. The
length of a smooth curve with respect to this metric is the integral of
$||\delta q||$ along the curve. The proof of the Theorem on Local Ergodicity
[S-Ch(1987)] shows that the $\delta q$-metric is the relevant notion of
distance on the local orthogonal manifolds $\Sigma$, also being in good
harmony with the tubular distance function $z_{tub}(\,.\,)$ defined above.

\medskip

The tangent vector $u(x_0)=(\delta\tilde q_0,\delta\tilde v_0)$ is defined
as follows:

For a large constant $L_0\gg 1$ (to be specified later), we select first a
non-collision time $\tilde c_3\gg 1$ in the following way: Thanks to our
hypothesis (5) in 3.1 and the hypersurface measure expansion theorem of
[Ch-Sim(2006)], the hypersurface measure of $S^t(J)$ grows at least linearly
in $t$, as $t\to\infty$. As a consequence, the distances between $J$ and
nearby points will shrink at least linearly in $t$, as $t\to\infty$. The
theorem of Appendix II below claims that for a large enough
(non-collision) time $\tilde c_3\gg 1$ the phase point $x_{\tilde
  c_3}=S^{\tilde c_3}x_0$ has a unit tangent vector $(\delta q_0,\delta
v_0)\in E^s(x_{\tilde c_3})$ such that the normalized tangent vector

$$
u(x_0)=
(\delta\tilde q_0,\delta\tilde v_0):=\frac{(DS^{-\tilde c_3})(\delta q_0,
\delta v_0)}{\left||(DS^{-\tilde c_3})(\delta q_0,\delta v_0)\right||}\in
E^s(x_{0}) \tag 3.8
$$
is transversal to $J$, and the expansion estimate

$$
\frac{||\delta\tilde q_{0}||}{||\delta\tilde q_{\tilde c_3}||}>2L_0
$$
holds true or, equivalently, we have the contraction estimate

$$
\frac{||\delta\tilde q_{\tilde c_3}||}{||\delta\tilde q_{0}||}<\frac{1}{2L_0},
\tag 3.9
$$
where 

$$
(\delta\tilde q_{\tilde c_3},\,\delta\tilde v_{\tilde c_3}):=
\left(DS^{\tilde c_3}\right)(\delta\tilde q_{0},\,\delta\tilde v_{0})=
\frac{(\delta q_0,\delta v_0)}{\left||(DS^{-\tilde c_3})
(\delta q_0,\delta v_0)\right||}.
$$

\medskip

\subheading{Remark} Almost every phase point $x_0$ of the hypersurface $J$
satisfies the hypotheses of Appendix II (on the connected collision
graphs). This is indeed so, since the proof of Theorem 6.1 of [Sim(1992)-I]
works without any essential change not only for singular phase points, but
also for the points of the considered exceptional manifold $J$. The only
important ingredient of that proof is the transversality of the spaces
$E^s(x)$ to $J$, provided by Lemma 3.3 above. According to that result,
typical phase points $x\in J$ (with respect to the hypersurface measure of
$J$) indeed enjoy the above property of having infinitely many consecutive,
connected collision graphs on their forward orbit $S^{(0,\infty)}x_0$.

\medskip

We choose the orientation of the unit normal field $n(x)$ ($x\in J$) of $J$ in
such a way that $\langle n(x_0),\,(\delta\tilde q_0,\,\delta\tilde
v_0)\rangle<0$, and define the one-sided tubular neighborhood $U_\delta$ of
radius $\delta>0$ as the set of all phase points $\gamma_x(s)$, where $x\in
J$, $0\le s<\delta$. Here $\gamma_x(\,.\,)$ is the geodesic line passing
through $x$ (at time zero) with the initial velocity $n(x)$, $x\in J$. The
radius (thickness) $\delta>0$ here is a variable, which will eventually tend
to zero. We are interested in getting useful asymptotic estimates for certain
subsets of $U_\delta$, as $\delta\to 0$.

Our main working domain will be the set

$$
\aligned
D_0=\Big\{&y\in U_{\delta_0}\setminus J\big|\; y\not\in\bigcup\Sb t>0\endSb
S^{-t}\left(\Cal S\Cal R^-\right),\;\exists\text{ a sequence} \\
&t_n\nearrow\infty\text{ such that }
S^{t_n}y\in U_{\delta_0}\setminus J, \quad n=1,2,\dots\Big\},
\endaligned
\tag 3.10
$$
a set of full $\mu$-measure in $U_{\delta_0}$. We will use the shorthand
notation $U_0=U_{\delta_0}$ for a fixed, small value $\delta_0$.

\medskip

On any manifold $\Sigma_0^t(y)\cap U_0$ we define the smooth
field $\Cal X_{y,t}(y')$ ($y'\in \Sigma_0^t(y)\cap U_0$) of unit tangent
vectors of $\Sigma_0^t(y)\cap U_0$ as follows:

$$
\Cal X_{y,t}(y')=\frac{\Pi_{y,t,y'}\left((\delta\tilde q_0,\,\delta
\tilde v_0)\right)}
{\left\Vert\Pi_{y,t,y'}\left((\delta\tilde q_0,\,\delta
\tilde v_0)\right)\right\Vert},
\tag 3.11
$$
where $\Pi_{y,t,y'}$ denotes the orthogonal projection of 
$\Bbb R^d\oplus\Bbb R^d$ onto the tangent space of $\Sigma_0^t(y)$ at the
point $y'\in\Sigma_0^t(y)\cap U_0$. 

\medskip

If, in the construction of the manifolds $\Sigma_0^t(y)$, the time $t$ is
large enough, i. e. $t\ge c_3$ for a suitably large constant 
$c_3\gg\tilde c_3$, the points $y$, $y'$ are close enough to $x_0$, and
$y'\in\Sigma_0^t(y)$, then the tangent space $\Cal T_{y'}\Sigma_0^t(y)$ will
be close enough to the tangent space $\Cal T_{x_0}\gamma^{s}(x_0)$ of the
local stable manifold $\gamma^{s}(x_0)$ of $x_0$, so that the projected copy
$\Cal X_{y,t}(y')$ of $(\delta\tilde q_0,\,\delta\tilde v_0)$ (featuring
(3.11)) will undergo a contraction by a factor of at least $L_0^{-1}$ between
time $0$ and $\tilde c_3$, let alone between time $0$ and $c_3$, that is,

$$
\frac{\left||DS^t\left(\Cal X_{y,t}(y')\right)\right||_q}
{\left||\left(\Cal X_{y,t}(y')\right)\right||_q}<L_0^{-1}
\tag 3.12
$$
for all $t\ge c_3$. We note that the tangent space $\Cal
T_{x_0}\gamma^{s}(x_0)$ of the local stable manifold $\gamma^{s}(x_0)$
makes sense, even if the latter object does not exist: this tangent space can
be obtained as the positive subspace of the operator $B(x_0)$ defined by the
continued fraction (2) in [S-Ch(1987)] or, equivalently, as the intersection
of the inverse images of stable cones of remote phase points on the forward
orbit of $x_0$.  All the necessary upper estimates for the mentioned angles
between the considered tangent spaces follow from the well known result
stating that the difference (in norm) between the second fundamental forms of
the $S^t$-images ($t>0$) of two local, convex orthogonal manifolds is at most
$1/t$, see, for instance, inequality (4) in [Ch(1982)]. These facts imply, in
particular, that the vector in the numerator of (3.11) is actually very close
to $(\delta\tilde q_0,\delta\tilde v_0)$, thus its magnitude is almost one.

\medskip

For any $y\in D_0$ let $t_k=t_k(y)$ ($0<t_1<t_2<\dots$) be the time of the
$k$-th collision $\sigma_k$ on the forward orbit $S^{[0,\infty)}y$ of $y$.
Assume that the time $t$ in the construction of $\Sigma_0^t(y)$
and $\Cal X_{y,t}$ is between $\sigma_{k-1}$ and $\sigma_k$, i. e. 
$t_{k-1}(y)<t<t_k(y)$. We define the smooth curve $\rho_{y,t}=\rho_{y,t}(s)$
(with the arc length parametrization $s$, $0\le s\le h(y,t)$) as the maximal
integral curve of the vector field $\Cal X_{y,t}$ emanating from $y$ and not
intersecting any forward singularity of order $\le k$, i. e.

$$
\cases
&\rho_{y,t}(0)=y, \\
&\frac{d}{ds}\rho_{y,t}(s)=\Cal X_{y,t}\left(\rho_{y,t}(s)\right), \\
&\rho_{y,t}(\,.\,)\text{ does not intersect any singularity of order }\le k, \\
&\rho_{y,t}\text{ is maximal among all curves with the above properties.}
\endcases
\tag 3.13
$$

We remind the reader that a phase point $x$ lies on a singularity of order $k$
($k\in\Bbb N$) if and only if the $k$-th collision on the forward orbit
$S^{(0,\infty)}x$ is a singular one. It is also worth noting here that, as it
immediately follows from definition (3.6) and (3.13), the curve $\rho_{y,t}$
can only terminate at a boundary point of the manifold $\Sigma_0^t(y)\cap
U_0$.

\medskip

\subheading{Remark 3.14} From now on, we will use the notations
$\Sigma_0^k(y)$, $\Cal X_{y,k}$, and $\rho_{y,k}$ for $\Sigma_0^{t_k^*}(y)$,
$\Cal X_{y,t_k^*}$, and $\rho_{y,t_k^*}$, respectively, where
$t_k^*=t_k^*(y)=\frac{1}{2}\left(t_{k-1}(y)+t_k(y)\right)$.

\medskip

Due to these circumstances, the curves $\rho_{y,t_k^*}=\rho_{y,k}$ can now
terminate at a point $z$ such that $z$ is not on any singularity of order at
most $k$ and $S^{t_k^*}z$ is a boundary point of $\Sigma_{t_k^*}^{t_k^*}(y)$,
so that at the point $S^{t_k^*}z$ the manifold $\Sigma_{t_k^*}^{t_k^*}(y)$
touches the boundary of the phase space in a nonsingular way. This means that,
when we continuously move the points $\rho_{y,k}(s)$ by varying the parameter
$s$ between $0$ and $h(y,k)$, either the time $t_k\left(\rho_{y,k}(s)\right)$
or the time $t_{k-1}\left(\rho_{y,k}(s)\right)$ becomes equal to
$t_k^*=t_k^*(y)$ when the parameter value $s$ reaches its maximum value
$h(y,k)$. The length of the curve $\rho_{y,k}$ is at most $\delta_0$, and an
elementary geometric argument shows that the time of collision
$t_k(\rho_{y,k}(s))$ (or $t_{k-1}(\rho_{y,k}(s))$) can only change by at most
the amount of $c^*\sqrt{\delta_0}$, as $s$ varies between $0$ and $h(y,k)$.
(Here $c^*$ is an absolute constant.) Thus, we get that the unpleasant
situation mentioned above can only occur when the difference
$t_k(y)-t_{k-1}(y)$ is at most $c^*\sqrt{\delta_0}$.  These collisions have to
be and will be excluded as stopping times $k_2(y)$, $\overline{t}_2(y)$ and
$\overline{k}_1(y)$ in the proof below. Still, everything works by the main
result of [B-F-K(1998)], which states that there is a large positive integer
$n_0$ and a small number $\beta>0$ such that amongst any collection of $n_0$
consecutive collisions there are always two neighboring ones separated from
each other (in time) by at least $\beta$. Taking $c^*\sqrt{\delta_0}<\beta$
shows that the badly behaved collisions -- described above -- can indeed
be excluded from our construction.

\medskip

As far as the terminal point $\rho_{y,k}(h(y,k))$ of $\rho_{y,k}$ is
concerned, there are exactly three, mutually exclusive possiblities for this 
point:

\medskip

(A) $\rho_{y,k}(h(y,k))\in J$ and this terminal point does not belong to any
forward singularity of order $\le k$,

\medskip

(B) $\rho_{y,k}(h(y,k))$ lies on a forward singularity of order $\le k$,

\medskip

(C) the terminal point $\rho_{y,k}(h(y,k))$ does not lie on any singularity
of order $\le k$ but lies on the part of the boundary
$\partial U_0$ of $U_0$ different from $J=J\times\{0\}$ and
$J\times\{\delta_0\}$.

\medskip

\subheading{Remark 3.15} Under the canonical identification $U_0\cong
J\times[0,\,\delta_0)$ of $U_0$ via the geodesic lines perpendicular to $J$,
the above mentioned part of $\partial U_0$ (the "side" of $U_0$) corresponds
to $\partial J\times[0,\,\delta_0)$. Therefore, the set of points with
property (C) inside a layer $U_\delta$ ($\delta\le\delta_0$) will have
$\mu$-measure $o(\delta)$ (actually, of order $\delta^2$), and this set will
be negligible in our asymptotic measure estimations, as $\delta\to 0$. The
reason why these sets are negligible, is that in the indirect proof of Main
Lemma 3.5, a contradiction will be obtained (at the end of \S3) by comparing
the measures of certain sets, whose measures are of order
$\text{const}\cdot\delta$.  That is why in the future we will not be dealing
with any phase point with property (C).

\medskip

Should (B) occur for some value of $k$ ($k\ge 2$), the minimum of all such
integers $k$ will be denoted by $\overline{k}=\overline{k}(y)$. The exact
order of the forward singularity on which the terminal point
$\rho_{y,\overline{k}}\left(h(y,\overline{k})\right)$ lies is denoted by
$\overline{k}_1=\overline{k}_1(y)$ ($\le\overline{k}(y)$).  If (B) does not
occur for any value of $k$, then we take
$\overline{k}(y)=\overline{k}_1(y)=\infty$.

We can assume that the manifold $J$ and its one-sided tubular neighborhood 
$U_0=U_{\delta_0}$ are already so small that for any $y\in U_0$ 
no singularity of $S^{(0,\infty)}y$ can take place at the first
collision, so the indices $\overline{k}$ and $\overline{k}_1$ above are
automatically at least $2$. For our purposes the important index will be
$\overline{k}_1=\overline{k}_1(y)$ for phase points $y\in D_0$. 

\medskip

\subheading{Remark 3.16. Refinement of the construction} Instead of selecting
a single contracting unit vector $(\delta\tilde q_0,\,\delta\tilde v_0)$ in
(3.8), we should do the following: Choose a compact set $K_0\subset B$ with
the property

$$
\frac{\mu_1(K_0)}{\mu_1(J)}>1-10^{-6}.
$$
Now the running point $x\in K_0$ will play the role of $x_0$ in the
construction of the contracting unit tangent vector
$u(x):=(\delta\tilde q_0,\,\delta\tilde v_0)\in E^s(x)$ on the
left-hand-side of (3.8). For every $x\in K_0$ there is a small, open ball
neighborhood $B(x)$ of $x$ and a big threshold $\tilde c_3(x)\gg 1$ such that 
the contraction estimate (3.9) and the transversality to $J$ 
hold true for $u(y)$ and $\tilde c_3=\tilde c_3(x)$ for all $y\in B(x)$.

Exactly the same way as earlier, one can also achieve that the weaker
contraction estimate $L_0^{-1}$ of (3.12) holds true not only for $t\ge c_3$
and $u(x)$, but also for any projected copy of it appearing in (3.11) and
(3.12), provided that $y,\, y'\in B(x)$, and $t\ge c_3(x)$.

Now select a finite subcover $\bigcup_{i=1}^n B(x_i)$ of $K_0$, and
replace $J$ by $J_1=J\cap\bigcup_{i=1}^n B(x_i)$, $U_\delta$ by
$U'_\delta=U_\delta\cap\bigcup_{i=1}^n B(x_i)$ (for
$\delta\le\delta_0$) and, finally, choose the threshold $c_3$ to be
the maximum of all thresholds $c_3(x_i)$ for $i=1,2,\dots,n$. In this
way the assertion of Corollary 3.18 below will be true.

We note that the new exceptional manifold $J_1$ is no longer so nicely "round
shaped" as $J$, but it is still pretty well shaped, being a
domain in $J$ with a piecewise smooth boundary.

The reason why we cannot switch completely to a round and much smaller
manifold $B(x)\cap J$ is that the measure $\mu_1(J)$ should be kept bounded
from below after having fixed $L_0$, see the requirement 3 in Appendix I.

\medskip

\subheading{Remark 3.17} When defining the returns of a forward orbit to
$U_\delta$, we used to say that ``before every new return the orbit must first
leave the set $U_\delta$''.  Since the newly obtained $J$ is no longer round
shaped as it used to be, the above phrase is not satisfactory any longer.
Instead, one should say that the orbit leaves even the $\kappa$-neighborhood
of $U_\delta$, where $\kappa$ is two times the diameter of the original $J$.
This guarantees that not only the new $U_\delta$, but also the original
$U_\delta$ will be left by the orbit, so we indeed are dealing with a genuine
return. This note also applies to two more shrinkings of $J$ that will
take place later in the proof.

\medskip

In addition, it should be noted that, when constructing the vector field in
(3.11) and the curves $\rho_{y,t}$, an appropriate directing vector
$u(x_i)$ needs to be chosen for (3.11). To be definite and not arbitrary, a
convenient choice is the first index $i\in\{1,2,\dots,n\}$ for which
$y\in B(x_i)$. In that way the whole curve $\rho_{y,t}$ will stay in the
slightly enlarged ball $B'(x_i)$ with double the radius of $B(x_i)$, and one
can organize things in such a way that the required contraction estimates of (3.12)
be still true, even in these enlarged balls.

In the future, a bit sloppily, $J_1$ will be denoted by $J$, and
$U'_\delta$ by $U_\delta$.

\medskip

As an immediate corollary of (3.12), the uniform transvarsality of the field
$\Cal X_{y,t}(y')$ to $J$ and Remark 3.16, we get

\medskip

\subheading{Corollary 3.18} For the given sets $J$, $U_0$, and the
large constant $L_0$ we can select the threshold $c_3>0$ large enough so
that for any point $y\in D_0$ any time $t$ with 
$c_3\le t<t_{\overline{k}_1(y)}(y)$ the $\delta q$-expansion rate of $S^t$
between the curves $\rho_{y,\overline{k}(y)}$ and
$S^t\left(\rho_{y,\overline{k}(y)}\right)$ is less than $L_0^{-1}$, i. e. for
any tangent vector $(\delta q_0,\,\delta v_0)$ of $\rho_{y,\overline{k}(y)}$ we
have 

$$
\frac{||\delta q_t||}{||\delta q_0||}<L_0^{-1},
$$
where $(\delta q_t,\,\delta v_t)=(DS^t)(\delta q_0,\,\delta v_0)$.

\medskip

\subheading{Remark 3.19}
The reason why there is no expansion from time $c_3$ until
time $t$ is that all the image curves
$S^{\tau}\left(\rho_{y,\overline{k}(y)}\right)$ ($c_3\le\tau\le t$) are
concave, according to the construction of the curve
$\rho_{y,\overline{k}(y)}$.

\medskip

An immediate consequence of the previous result is

\medskip

\subheading{Corollary 3.20} For any $y\in D_0$ with $\overline{k}(y)<\infty$
and $t_{\overline{k}_1(y)-1}(y)\ge c_3$, and for any $t$ with
$t_{\overline{k}_1(y)-1}(y)<t<t_{\overline{k}_1(y)}(y)$, we have 

$$
z_{tub}\left(S^ty\right)<L_0^{-1}l_q\left(\rho_{y,\overline{k}(y)}\right)
<\frac{c_4}{L_0}\text{dist}(y,\,J),
\tag 3.21
$$
where $l_q\left(\rho_{y,\overline{k}(y)}\right)$ denotes the $\delta q$-length
of the curve $\rho_{y,\overline{k}(y)}$, and $c_4>0$ is a constant,
independent of $L_0$ or $c_3$, depending only on the (asymptotic) angles 
between the curves $\rho_{y,\overline{k}(y)}$ and $J$. 

\medskip

\subheading{Proof} The manifold $J$ and the curves $\rho_{y,\overline{k}(y)}$
are uniformly transversal, as it follows immediately from the
uniformity of the transversality of the field $\Cal X_{y,t}(y')$ to $J$.  This
is why the above constant $c_4$, independently of $L_0$, exists. \qed

\medskip

By further shrinking the exceptional manifold $J$ a little bit,
and by selecting a suitably
thin, one-sided neighborhood $U_1=U_{\delta_1}$ of $J$, we can achieve that
the open $2\delta_1$-neighborhood of $U_1$ (on the same side of $J$ as
$U_0$ and $U_1$) is a subset of $U_0$.

For a varying $\delta$, $0<\delta\le\delta_1$, we introduce the layer

$$
\aligned
\overline{U}_\delta=\Big\{&y\in(U_\delta \setminus U_{\delta/2})\cap
D_0\big|\;\exists\text{ a sequence }t_n\nearrow\infty \\
&\text{such that } S^{t_n}y\in(U_\delta \setminus U_{\delta/2})
\text{ for all }n\Big\}.
\endaligned
\tag 3.22
$$
Since almost every point of the layer 
$(U_\delta \setminus U_{\delta/2})\cap D_0$ returns infinitely often to
this set and the asymptotic equation 

$$
\mu\left((U_\delta \setminus U_{\delta/2})\cap D_0\right)\sim
\frac{\delta}{2}\mu_1(J)
\tag 3.23
$$
holds true, we get the asymptotic equation 

$$
\mu\left(\overline{U}_\delta\right)
\sim\frac{\delta}{2}\mu_1(J).
\tag 3.24
$$

We will need the following subsets of $\overline{U}_\delta$:

$$
\aligned
\overline{U}_\delta(c_3)&=\left\{y\in\overline{U}_\delta\big|\;
t_{\overline{k}_1(y)-1}(y)\ge c_3\right\}, \\
\overline{U}_\delta(\infty)&=\left\{y\in\overline{U}_\delta\big|\;
\overline{k}_1(y)=\infty\right\}.
\endaligned
\tag 3.25
$$
Here $c_3$ is the constant from Corollary 3.18, the exact value of which will
be specified later, at the end of the proof of Main Lemma 3.5. Note that in
the first line of (3.25) the case $\overline{k}_1(y)=\infty$ is included.  By
selecting the pair of sets $(U_1,\,J)$ small enough, we can assume that

$$
z_{tub}(y)>c_4\delta_1 \quad\forall y\in U_1.
\tag 3.26
$$
This inequality guarantees that the collision time $t_{\overline{k}_1(y)}(y)$
($y\in\overline{U}_\delta$) cannot be near any return time of $y$ to the
layer $(U_\delta\setminus U_{\delta/2})$, for $\delta\le\delta_1$,
provided that $y\in\overline{U}_\delta(c_3)$.
More precisely, the whole orbit segment 
$S^{[-\tau(-z),\,\tau(z)]}z$ will be disjoint from $U_1$, where
$z=S^ty$, $t_{\overline{k}_1(y)-1}(y)<t<t_{\overline{k}_1(y)}(y)$.

\medskip

Let us consider now the points $y$ of the set $\overline{U}_\delta(\infty)$.
We observe that for any point $y\in\overline{U}_\delta(\infty)$ the curves
$\rho_{y,k}(s)$ ($0\le s\le h(y,k)$) have a $C^2$-limiting curve
$\rho_{y,\infty}(s)$ ($0\le s\le h(y,\infty)$), with $h(y,\,k)\to
h(y,\,\infty)$, as $k\to\infty$.

Indeed, besides the concave, local orthogonal manifolds
$\Sigma_0^k(y)=\Sigma_0^{t_k^*}(y)$ of (3.6) (where
$t_k^*=t_k^*(y)=\frac{1}{2}\left(t_{k-1}(y)+t_k(y)\right)$), let us also
consider another type of concave, local orthogonal manifolds defined by
the formula

$$
\tilde\Sigma_0^k(y)=
\tilde\Sigma_0^{t_k^*}(y)=SC_y\left(S^{-t_k^*}\left(SC_{y_{t_k^*}}\left\{y'\in
\bold M\big|\; q(y')=q(y_{t_k^*})\right\}\right)\right),
$$
the so called "candle manifolds", containing the phase point
$y\in\overline{U}_\delta(\infty)$ in their interior. It was proved in \S3
of [Ch(1982)] that the second fundamental forms
$B\left(\Sigma_0^k(y),\,y\right)\le 0$ are monotone
non-increasing in $k$, while the second fundamental forms
$B\left(\tilde\Sigma_0^k(y),\,y\right)<0$ are monotone
increasing in $k$, so that

$$
B\left(\tilde\Sigma_0^k(y),\,y\right)<B\left(\Sigma_0^k(y),\,y\right)
$$
is always true. It is also proved in \S3 of [Ch(1982)] that

$$
\lim_{t\to\infty} B\left(\tilde\Sigma_0^k(y),\,y\right)=
\lim_{t\to\infty} B\left(\Sigma_0^k(y),\,y\right):=B_\infty(y)<0
$$
uniformly in $y$, and these two-sided, monotone curvature limits give rise to
uniform $C^2$-convergences

$$
\lim_{t\to\infty}\Sigma_0^k(y)=\Sigma_0^\infty(y),\quad
\lim_{t\to\infty}\tilde\Sigma_0^k(y)=\Sigma_0^\infty(y),
$$
and the limiting manifold $\Sigma_0^\infty(y)$ is the local stable invariant
manifold $\gamma^{s}(y)$ of $y$, once it contains $y$ in its smooth
part. These monotone, two-sided limit relations, together with the definition
of the curves $\rho_{y,{t_k^*}}=\rho_{y,k}$ prove the
existence of the $C^2$-limiting curve
$\rho_{y,\infty}=\lim_{k\to\infty}\rho_{y,k}$,\;
$h(y,k)\to h(y,\infty)$, as $k\to\infty$. They also prove the inclusion
$\rho_{y,\infty}\left([0,\,h(y,\infty)]\right)\subset\gamma^{s}(y)$.

\medskip

\subheading{Lemma 3.27} 
$\mu\left(\overline{U}_\delta\setminus\overline{U}_\delta(c_3)\right)
=o(\delta)$, as $\delta\to 0$.

\medskip

\subheading{Proof} The points $y$ of the set 
$\overline{U}_\delta\setminus\overline{U}_\delta(c_3)$ have the property
$t_{\overline{k}_1(y)-1}(y)<c_3$. By doing another slight shrinking to $J$,
the same way as in Remark 3.16, we can achieve that 
$t_{\overline{k}_1(y)}(y)<2c_3$ for all 
$y\in\overline{U}_\delta\setminus\overline{U}_\delta(c_3)$, 
$0<\delta\le\delta_1$. This means that the terminal point
$\Pi(y)=\rho_{y,\overline{k}_1(y)}\left(h(y,k)\right)$ of the curve
$\rho_{y,\overline{k}_1(y)}$ lies on the singularity set

$$
\bigcup\Sb 0\le t\le 2c_3\endSb S^{-t}\left(\Cal S\Cal R^-\right),
$$
hence all points of the set
$\overline{U}_\delta\setminus\overline{U}_\delta(c_3)$ are at most at
the distance of $\delta$ from the singularity set mentioned above.

This singularity set is a compact collection of codimension-one, smooth
submanifolds (with boundaries), each of which is uniformly transversal to the
manifold $J$. This uniform transversality follows from Lemma 3.3 above, and
from the fact that the inverse images $S^{-t}(\Cal S\Cal R^-)$ ($t>0$) of
singularities can be smoothly foliated with local, concave orthogonal
manifolds. Thus, the $\delta$-neighborhood of this singularity set inside 
$\overline{U}_\delta$ clearly has $\mu$-measure $o(\delta)$,
actually, of order $\le\text{const}\cdot\delta^2$. \qed

\medskip

For any point $y\in\overline{U}_\delta(\infty)$ we define the return time
$\overline{t}_2=\overline{t}_2(y)$ as the infimum of all numbers $t_2>c_3$
for which there exists another number $t_1$, $0<t_1<t_2$, such that 
$S^{t_1}y\not\in \tilde{U}_0$ and 
$S^{t_2}(y)\in\left(U_\delta\setminus U_{\delta/2}\right)\cap D_0$. Let 
$k_2=k_2(y)$ be the unique natural number for which 
$t_{k_2-1}(y)<\overline{t}_2(y)<t_{k_2}(y)$.

\medskip

\subheading{Lemma 3.28} For any point $y\in\overline{U}_\delta(\infty)$
the projection 

$$
\Pi(y):=\rho_{y,k_2(y)}\left(h(y,k_2(y))\right)
$$
is a forward singular point of $J$.

\medskip

\subheading{Proof} Assume that the forward orbit of $\Pi(y)$ is non-singular.
The distance $\text{dist}(S^{\overline{t}_2}y,\,J)$ between
$S^{\overline{t}_2}y$ and $J$ is bigger than $\delta/2$.  According to the
contraction result 3.18, if the contraction factor $L_0^{-1}$ is chosen small
enough, the distance between $S^{\overline{t}_2}\left(\Pi(y)\right)$ and $J$
stays bigger than $\delta/4$, so $S^{\overline{t}_2}\left(\Pi(y)\right)\in
U_0\setminus J$ will be true.  This means, on the other hand, that the forward orbit of
$\Pi(y)$ is sufficient, according to (4)/(ii) of Definition 3.1. However, this
is impossible, due to our standing assumption $A=\emptyset$. \qed

\medskip

\subheading{Lemma 3.29} The set $\overline{U}_\delta(\infty)$ is actually empty.

\medskip

\subheading{Proof} Just observe that in the previous proof the whole curve
$\rho_{y,k_2(y)}$ can be slightly perturbed (in the $C^\infty$ topology, for
example), so that the perturbed curve $\tilde{\rho}_y$ emanates from $y$ and
terminates on a non-singular point $\tilde{\Pi}(y)$ of $J$ (near $\Pi(y)$),
so that the curve $\tilde{\rho}_y$ still "lifts" the point $\tilde{\Pi}(y)$
up to the set $\left(U_\delta\setminus U_{\delta/4}\right)\cap D_0$ if we
apply $S^{\overline{t}_2}$. This proves the existence of a non-singular,
sufficient phase point $\tilde{\Pi}(y)\in A$, which is impossible by our
standing assumption $A=\emptyset$. Hence 
$\overline{U}_\delta(\infty)=\emptyset$. \qed

\medskip

Next we need a useful upper estimate for the $\mu$-measure of the set
$\overline{U}_\delta(c_3)$ as $\delta\to 0$. We will classify the points
$y\in\overline{U}_\delta(c_3)$ according to whether $S^ty$ returns to the
layer $\left(U_\delta\setminus U_{\delta/2}\right)\cap D_0$ (after first
leaving it, of course) before the time $t_{\overline{k}_1(y)-1}(y)$ or not.
Thus, we define the sets

$$
\aligned
E_\delta(c_3)=&\big\{y\in\overline{U}_\delta(c_3)
\big|\; \exists\; 0<t_1<t_2<t_{\overline{k}_1(y)-1}(y) \\
&\text{such that }S^{t_1}y\not\in\tilde{U}_0,\; S^{t_2}y\in
\left(U_\delta\setminus U_{\delta/2}\right)\cap D_0\big\}, \\
F_\delta(c_3)=&\overline{U}_\delta(c_3)\setminus E_\delta(c_3).
\endaligned
\tag 3.30
$$
Recall that the threshold $t_{\overline{k}_1(y)-1}(y)$, being a collision
time, is far from any possible return time $t_2$ to the layer
$\left(U_\delta\setminus U_{\delta/2}\right)\cap D_0$, see the remark
right after (3.26).

Now we will be doing the "slight shrinking" trick of Remark 3.16 the third
(and last) time. We slightly further decrease $J$ to obtain a smaller 
$J_1$ with almost the same $\mu_1$-measure. Indeed, by using property (6)
of 3.1, inside the set $J\cap B$ we choose a compact set $K_1$ for which

$$
\frac{\mu_1(K_1)}{\mu_1(J)}>1-10^{-6},
$$
and no point of $K_1$ ever returns to $J$. For each point $x\in K_1$ the
distance between the orbit segment $S^{[a_0, c_3]}x$ and $J$ is at least
$\epsilon(x)>0$. Here $a_0$ is needed to guarantee that we certainly drop
the initial part of the orbit, which still stays near $J$, and $c_3$ was chosen
earlier. By the non-singularity of the orbit segment $S^{[a_0, c_3]}x$ and
by continuity, the point $x\in K_1$ has an open ball neighborhood $B(x)$
of radius $r(x)>0$ such that for every $y\in B(x)$ the orbit segment
$S^{[a_0, c_3]}y$ is non-singular and stays away from $J$ by at least
$\epsilon(x)/2$. Choose a finite covering
$\bigcup_{i=1}^n B(x_i)\supset K_1$ of $K_1$, replace $J$ and $U_\delta$ by
their intersections with the above union (the same way as it was done in
Remark 3.16), and fix the threshold value of $\delta_1$ so that

$$
\delta_1<\frac{1}{2}\min\{\epsilon(x_i)|\, i=1,2,\dots,n\}.
$$
In the future we again keep the old notations $J$ and $U_\delta$ for these 
intersections. In this way we achieve that the following statement be true:

$$
\cases &\text{any return time }t_2\text{ of any point }y\in
\left(U_\delta\setminus U_{\delta/2}\right)\cap D_0\text{ to} \\
&\left(U_\delta\setminus U_{\delta/2}\right)\cap D_0 \text{ is always greater
  than }c_3\text{ for } 0<\delta\le\delta_1.
\endcases
\tag 3.31
$$
Just as in the paragraph before Lemma 3.28,
for any phase point $y\in E_\delta(c_3)$ we define the return time 
$\overline{t}_2=\overline{t}_2(y)$ as the infimum of all the return times
$t_2$ of $y$ featuring (3.30). By using this definition of
$\overline{t}_2(y)$, formulas (3.30)--(3.31), and the contraction result
3.18, we easily get

\medskip

\subheading{Lemma 3.32} If the contraction coefficient $L_0^{-1}$ in 3.18 is
chosen suitably small, for any point $y\in E_\delta(c_3)$ the
projected point

$$
\Pi(y):=\rho_{y,\overline{t}_2(y)}\left(h(y,\overline{t}_2(y))\right)\in J
\tag 3.33
$$
is a forward singular point of $J$.

\medskip

\subheading{Proof} Since $\overline{t}_2(y)<t_{\overline{k}_1(y)-1}(y)$, we
get that, indeed, $\Pi(y)\in J$. Assume that the forward orbit of $\Pi(y)$
is non-singular.

Since $S^{\overline{t}_2(y)}y\in\overline{\left(U_\delta\setminus 
U_{\delta/2}\right)\cap D_0}$, we obtain that 
$\text{dist}\left(S^{\overline{t}_2(y)}y,\, J\right)\ge\delta/2$.
On the other hand, by using (3.31) and Corollary 3.18,
we get that for a small enough
contraction coefficient $L_0^{-1}$ the distance between 
$S^{\overline{t}_2(y)}y$ and $S^{\overline{t}_2(y)}\left(\Pi(y)\right)$
is less than $\delta/4$. (The argument is the same as in the proof of 
Lemma 3.28.) In this way we obtain that 
$S^{\overline{t}_2(y)}\left(\Pi(y)\right)\in U_0\setminus J$, so 
$\Pi(y)\in A$, according to condition (4)/(ii) in 3.1, thus contradicting to
our standing assumption $A=\emptyset$. This proves that, indeed, $\Pi(y)$
is a forward singular point of $J$. \qed

\medskip

\subheading{Lemma 3.34} The set $E_\delta(c_3)$ is actually empty.

\medskip

\subheading{Proof} The proof will be analogous with the proof of Lemma 3.29
above. Indeed, we observe that in the previous proof for any point 
$y\in E_\delta(c_3)$ the curve $\rho_{y,\overline{t}_2(y)}$ can be slightly
perturbed (in the $C^\infty$ topology), so that the perturbed curve
$\tilde{\rho}_y$ emanates from $y$ and terminates on a non-singular point
$\tilde{\Pi}(y)$ of $J$, so that the curve $\tilde{\rho}_y$
still "lifts" the point $\tilde{\Pi}(y)$ up to the set
$\left(U_\delta\setminus U_{\delta/2}\right)\cap D_0$ if we apply
$S^{\overline{t}_2}$. This means, however, that the terminal point 
$\tilde{\Pi}(y)$ of $\tilde{\rho}_y$ is an element of the set $A$, violating
our standing assumption $A=\emptyset$. This proves that no point
$y\in E_\delta(c_3)$ exists. \qed

\medskip

For the points $y\in F_\delta(c_3)=\overline{U}_\delta(c_3)$ we define the
projection $\tilde{\Pi}(y)$ by the formula

$$
\tilde{\Pi}(y):=S^{t_{\overline{k}_1(y)-1}(y)}y\in\partial\bold M.
\tag 3.35
$$
Now we prove

\medskip

\subheading{Lemma 3.36} For the measure 
$\nu\left(\tilde{\Pi}\left(F_\delta(c_3)\right)\right)$ of the projected set
$\tilde{\Pi}\left(F_\delta(c_3)\right)\subset\partial\bold M$ we have the upper
estimate 

$$
\nu\left(\tilde{\Pi}\left(F_\delta(c_3)\right)\right)\le c_2c_4L_0^{-1}\delta,
$$
where $c_2>0$ is the geometric constant (also denoted by $c_2$) in Lemma 2 of 
[S-Ch(1987)] or in Lemma 4.10 of [K-S-Sz(1990)-I], $c_4$ is the constant in
(3.21) above, and $\nu$ is the natural $T$-invariant measure on 
$\partial\bold M$ that can be obtained by projecting the Liouville measure
$\mu$ onto $\partial\bold M$ along the billiard flow.

\medskip

\subheading{Proof} Let $y\in F_\delta(c_3)$. From the inequality 
$t_{\overline{k}_1(y)-1}(y)\ge c_3$ and from Corollary 3.20 we conclude that 
$z_{tub}\left(\tilde{\Pi}(y)\right)<c_4L_0^{-1}\delta$. This inequality, along with
the fundamental measure estimate of Lemma 2 of [S-Ch(1987)] (see also Lemma
4.10 in [K-S-Sz(1990)-I]) yield the required upper estimate for
$\nu\left(\tilde{\Pi}\left(F_\delta(c_3)\right)\right)$. \qed

\medskip

The next lemma claims that the projection $\tilde{\Pi}:\;
F_\delta(c_3)\to\partial\bold M$ (considered here only on the set
$F_\delta(c_3)=\overline{U}_\delta(c_3)$) is "essentially one-to-one", from
the point of view of the Poincar\'e section.

\medskip

\subheading{Lemma 3.37} Suppose that $y_1,\, y_2\in F_\delta(c_3)$ are
non-periodic points ($\delta\le\delta_1$),
and $\Pi(y_1)=\Pi(y_2)$. We claim that $y_1$ and $y_2$
belong to an orbit segment $S$ of the billiard flow lying entirely in the
one-sided neighborhood $\tilde{U}_0$ of $J$ and, consequently, the
length of the segment $S$ is at most $1.1\text{diam}(J)$.

\medskip

\subheading{Remark} We note that in the length estimate 
$1.1\text{diam}(J)$ above, the coefficient $1.1$ could be replaced by any
number bigger than $1$, provided that the parameter $\delta>0$ is small
enough.

\medskip

\subheading{Proof} The relation $\Pi(y_1)=\Pi(y_2)$ implies that $y_1$ and
$y_2$ belong to the same orbit, so we can assume, for example, that
$y_2=S^ay_1$ with some $a>0$. We need to prove that $S^{[0,a]}y_1\subset\tilde U_0$.
Assume the opposite, i. e. that there is a number $t_1$, $0<t_1<a$, such that
$S^{t_1}y_1\not\in\tilde{U}_0$. This, and the relation
$S^ay_1\in\left(U_\delta\setminus U_{\delta/2}\right)\cap D_0$ mean that the
first return of $y_1$ to $\left(U_\delta\setminus U_{\delta/2}\right)\cap D_0$
occurs not later than at time $t=a$. On the other hand, since
$\Pi(y_1)=\Pi(S^ay_1)$ and $y_1$ is non-periodic, we get that
$t_{\overline{k}_1(y_1)-1}(y_1)>a$, see (3.35). The obtained inequality
$t_{\overline{k}_1(y_1)-1}(y_1)>a\ge\overline{t}_2(y_1)$, however, contradicts
to the definition of the set $F_\delta(c_3)$, to which $y_1$ belongs as an
element, see (3.30). The upper estimate $1.1\text{diam}(J)$ for the length of
$S$ is an immediate corollary of the containment $S\subset\tilde U_0$. \qed

\medskip

As a direct consequence of lemmas 3.36 and 3.37, we obtain

\medskip

\subheading{Corollary 3.38} For all small enough $\delta>0$, the inequality

$$
\mu\left(F_\delta(c_3)\right)\le 1.1c_2c_4L_0^{-1}\delta\text{diam}(J)
$$
holds true.

\bigskip

\subheading{Finishing the Indirect Proof of Main Lemma 3.5}

\medskip

It follows immediately from Lemma 3.27 and corollaries 3.29, 3.34, and 3.38
that 

$$
\mu\left(\overline{U}_\delta\right)\le 1.2c_2c_4\text{diam}(J)L_0^{-1}\delta
$$
for all small enough $\delta>0$. This fact, however, contradicts to (3.24)
if $L_0^{-1}$ is selected so small that

$$
1.2c_2c_4\text{diam}(J)L_0^{-1}<\frac{1}{4}\mu_1(J^*),
$$
where $J^*$ stands for the original exceptional manifold before the three
slight shrinkings in the style of Remark 3.16. Clearly,
$\mu_1(J)>(1-10^{-5})\mu_1(J^*)$. The obtained contradiction finishes the
indirect proof of Main Lemma 3.5. \qed

\bigskip \bigskip

\heading
\S4. Proof of Ergodicity \\
Induction on $N$
\endheading

\bigskip

By using several results of Sinai [Sin(1970)], Chernov-Sinai [S-Ch(1987)], and
Kr\'amli-Sim\'anyi-Sz\'asz, in this section we finally prove the ergodicity
(hence also the Bernoulli property; see Sinai's results in [Sin(1968)],
[Sin(1970)], and [Sin(1979)] for the K-mixing property, and then by
Chernov-Haskell [C-H(1996)] and Ornstein-Weiss [O-W(1998)] the Bernoulli
property follows from mixing) for every hard ball system $(\bold
M,\,\{S^t\},\,\mu)$, by carrying out an induction on the number $N$ ($\ge 2$)
of interacting balls.

\medskip

{\it In order for the proof to work, from now on we assume that the
Chernov-Sinai Ansatz is true for every hard ball system.}

\medskip

The base of the induction (i. e. the ergodicity of any two-ball system on a
flat torus) was proved in [Sin(1970)] and [S-Ch(1987)].

Assume now that $(\bold M,\,\{S^t\},\,\mu)$ is a given system of $N$ ($\ge 3$)
hard spheres with masses $m_1,m_2,\dots,m_N$ and radius $r>0$ on the flat unit
torus $\Bbb T^\nu=\Bbb R^\nu/\Bbb Z^\nu$ ($\nu\ge2$), as defined in \S2.
Assume further that the ergodicity of every such system is already proved to
be true for any number of balls $N'$ with $2\le N'<N$. We will carry out
the induction step by following the strategy for the proof laid down by Sinai
in [Sin(1979)] and polished in the series of papers [K-S-Sz(1989)],
[K-S-Sz(1990)-I], [K-S-Sz(1991)], and [K-S-Sz(1992)].

\medskip

By using the induction hypothesis, Theorem 5.1 of [Sim(1992)-I], together with
the slimness of the set $\Delta_2$ of doubly singular phase points, shows that
there exists a slim subset $S_1\subset\bold M$ of the phase space such that
for every $x\in\bold M\setminus S_1$ the point $x$ has at most one singularity
on its entire orbit $S^{(-\infty,\infty)}x$, and each branch of 
$S^{(-\infty,\infty)}x$ is not eventually splitting in any of the time
directions. By Corollary 3.26 and Lemma 4.2 of [Sim(2002)] there exists a
locally finite (hence countable) family of codimension-one, smooth,
exceptional submanifolds $J_i\subset\bold M$ such that for every point
$x\not\in\left(\bigcup_{i}J_i\right)\cup S_1$ the orbit of $x$ is sufficient
(geometrically hyperbolic). This means, in particular, that the considered
hard ball system $\flow$ is fully hyperbolic. 

By our standing assumption the Chernov-Sinai Ansatz (a global hypothesis
of the Theorem on Local Ergodicity by Chernov and Sinai, Theorem 5 in
[S-Ch(1987)], see also Corollary 3.12 in [K-S-Sz(1990)-I] and the main result
of [B-Ch-Sz-T(2002)]) is true, therefore, by the Theorem on Local Ergodicity,
an open neighborhood $U_x\ni x$ of any phase point
$x\not\in\left(\bigcup_{i}J_i\right)\cup S_1$ belongs to a single ergodic
component of the billiard flow. (Modulo the zero sets, of course.) Therefore,
the billiard flow $\{S^t\}$ has at most countably many, open ergodic
components $C_1,\,C_2,\,\dots$.

\medskip

\subheading{Remark} Note that theorem 5.1 of [Sim(1992)-I] 
(used above) requires the induction hypothesis as an assumption.

\medskip

Assume that, contrary to the statement of our theorem, the number of
ergodic components $C_1,\,C_2,\,\dots$ is more than one. The above
argument shows that, in this case, there exists a codimension-one,
smooth (actually analytic) submanifold $J\subset\bold
M\setminus\partial\bold M$ separating two different ergodic components
$C_1$ and $C_2$, lying on the two sides of $J$. By the Theorem on
Local Ergodicity for semi-dispersive billiards, no point of $J$ can
have a sufficient orbit. (Recall that sufficiency is clearly an open
property, so the existence of a sufficient point $y\in J$ would imply
the existence of a sufficient point $y'\in J$ with a non-singular
orbit.) By shrinking $J$, if necessary, we can achieve that the
infinitesimal Lyapunov function $Q(n)$ be separated from zero on $J$,
where $n$ is a unit normal field of $J$. By replacing $J$ with its
time-reversed copy

$$
-J=\left\{(q,v)\in\bold M\big|\; (q,-v)\in J\right\},
$$
if necessary, we can always achieve that $Q(n)\le -c_1<0$ uniformly across
$J$.

\medskip

There could be, however, a little difficulty in achieving the inequality
$Q(n)<0$ across $J$. Namely, it may happen that $Q(n_t)=0$ for every $t\in\Bbb
R$.  According to (7.2) of [Sim(2003)], the equation $Q(n_t)=0$ ($\forall\,
t\in\Bbb R$) implies that for the normal vector $n_t$ of $S^t(J)$ at
$x_t=S^tx_0$ one has $n_t=(0,\, w_t)$ for all $t\in\Bbb R$ and,
moreover, in the view of (7.5) of [Sim(2003)], $w_t^+=Rw_t^-$ is the
transformation law at any collision $x_t=(q_t,\, v_t)\in\partial\bold
M$. Furthermore, at every collision $x_t=(q_t,\, v_t)\in\partial\bold M$ the
projected tangent vector $V_1Rw_t^-=V_1w_t^+$ lies in the null space of the
operator $K$ (see also (7.5) in [Sim(2003)]), and this means that $w_0$ is a
neutral vector for the entire trajectory $S^{\Bbb R}y$, i. e. $w_0\in\Cal
N_0\left(S^{\Bbb R}y\right)$. On the other hand, this is impossible for the
following reason: any tangent vector $(\delta q,\delta v)$ from the space
$\Cal N_0\left(S^{\Bbb R}y\right)\times\Cal N_0\left(S^{\Bbb R}y\right)$ is
automatically tangent to the exceptional manifold $J$ (as a direct inspection
shows), thus for any normal vector $n=(z,w)\in\Cal T_x\bold M$ of a separating
manifold $J$ one has

$$
(z,\, w)\in\Cal N_0\left(S^{\Bbb R}y\right)^\perp\times\Cal N_0\left(
S^{\Bbb R}y\right)^\perp.
$$
The membership in this formula is, however, impossible with a nonzero vector
$w\in\Cal N_0\left(S^{\Bbb R}y\right)$.

\medskip

To make sure that the submanifold $J$ is neatly shaped (i. e. it fulfills (2)
of 3.1) is an obvious task. Condition (3) of 3.1 clearly holds true. We can
achieve (4) as follows: Select a base point $x_0\in J$ with a non-singular and
not eventually splitting forward orbit $S^{(0,\infty)}x_0$. This can be done
according to the transversality result 3.3 above (see also 7.12 in
[Sim(2003)]), and by using the fact that the points with an eventually
splitting forward orbit form a slim set in $\bold M$ (Theorem 5.1 of
[Sim(1992)-I]), henceforth a set of first category in $J$. After this, choose
a large enough time $T>0$ so that $S^Tx_0\not\in\partial\bold M$, and the
symbolic collision sequence $\Sigma_0=\Sigma\left(S^{[0,T]}x_0\right)$ is
combinatorially rich in the sense of Definition 3.28 of [Sim(2002)]. By
further shrinking $J$, if necessary, we can assume that
$S^T(J)\cap\partial\bold M=\emptyset$ and $S^T$ is smooth on $J$. Choose a
thin, tubular neighborhood $\tilde{U}_0$ of $J$ in $\bold M$ in such a way
that $S^T$ be still smooth across $\tilde{U}_0$, and define the set

$$
NS\left(\tilde{U}_0,\,\Sigma_0\right)=\left\{x\in\tilde{U}_0\big|\;\text{dim}
\Cal N_0\left(S^{[0,T]}x\right)>1\right\}
\tag 4.1
$$
of not $\Sigma_0$-sufficient phase points in $\tilde{U}_0$.  Clearly,
$NS\left(\tilde{U}_0,\,\Sigma_0\right)$ is a closed, algebraic set
containing $J$. We can assume that the selected (generic) base point
$x_0\in J$ belongs to the smooth part of the closed algebraic set
$NS\left(\tilde{U}_0,\,\Sigma_0\right)$.  This guarantees that
actually $J=NS\left(\tilde{U}_0,\,\Sigma_0\right)$, as long as the
manifold $J$ and its tubular neighborhood $\tilde{U}_0$ are selected
small enough, thus achieving property (4) of 3.1.

\medskip

\subheading{Proof of why property (6) of Definition 3.1 can be assumed}

\medskip

We recall that $J$ is a codimension-one, smooth manifold of non-sufficient 
phase points separating two open ergodic components, as described in (0)--(3)
at the end of \S3 of [Sim(2003)].

Let $P$ be the subset of $J$ containing all points with non-singular forward
orbit and recurring to $J$ infinitely many times.

\medskip

\subheading{Lemma 4.2} $\mu_1(P)=0$.

\medskip

\subheading{Proof} Assume that $\mu_1(P)>0$. Take a suitable Poincar\'e
section to make the time discrete, and consider the on-to-one first return map
$T:\;P\to P$ of $P$.  According to the measure expansion theorem for
hypersurfaces $J$ (with negative infinitesimal Lyapunov function $Q(n)$ for
their normal field $n$), proved in [Ch-Sim(2006)], the measure
$\mu_1\left(T(P)\right)$ is strictly larger than $\mu_1(P)$, though
$T(P)\subset P$. The obtained contradiction proves the lemma. \qed

\medskip

Next, we claim that the above lemma is enough for our purposes to prove (6) of
3.1. Indeed, the set $W\subset J$ consisting of all points $x\in J$ never
again returning to $J$ (after leaving it first, of course) has positive
$\mu_1$-measure by Lemma 4.2. Select a Lebesgue density base point $x_0\in W$
for $W$ with a non-singular forward orbit, and shrink $J$ at the very
beginning to such a small size around $x_0$ that the relative measure of $W$
in $J$ be bigger than $1-10^{-8}$.

\medskip

Finally, Main Lemma 3.5 asserts that $A\ne\emptyset$,
contradicting to our earlier statement that no point of $J$ is
sufficient. The obtained contradiction completes the inductive step of the
proof of the Theorem. \qed

\bigskip

\heading
Appendix I. The Constants of \S3
\endheading

\bigskip

In order to make the reading of \S3 easier, here we briefly describe
the hierarchy of the constants used there.

\medskip

1. The geometric constant $-c_1<0$ provides an upper estimate for the
   infinitesimal Lyapunov function $Q(n)$ of $J$ in (5) of Definition 3.1. It
   cannot be freely chosen in the proof of Main Lemma 3.5.

\medskip

2. The constant $c_2>0$ is present in the upper measure estimate of Lemma 2
   of [S-Ch(1987)], or Lemma 4.10 in [K-S-Sz(1990)-I]. It cannot be changed in
   the course of the proof of Main Lemma 3.5.

\medskip

3. The contraction coefficient $0<L_0^{-1}\ll 1$ plays a role all over \S3. It
must be chosen suitably small by selecting the time threshold $c_3\gg 1$ large
enough (see Corollary 3.20), after having fixed $\tilde U_0$, $\delta_0$, and $J$.
The phrase "suitably small" for $L_0^{-1}$ means that the inequality

$$
L_0^{-1}<\frac{0.25\mu_1(J^*)}{1.2c_2c_4\text{diam}(J)}
$$
should be true, see the end of \S3.

\medskip

4. The geometric constant $c_4>0$ of (3.21) bridges the gap between two
distances: the distance $\text{dist}(y,J)$ between a point $y\in D_0$ and
$J$, and the arc length $l_q\left(\rho_{y,\overline{k}(y)}\right)$.  It cannot
be freely chosen during the proof of Main Lemma 3.5.

\bigskip \bigskip

\heading
Appendix II \\
Expansion and Contraction Estimates 
\endheading

\bigskip \bigskip

For any phase point $x\in\bold M\setminus\partial\bold M$ with a
non-singular forward orbit $S^{(0,\infty)}x$ (and with at least one collision,
hence infinitely many collisions on it) we define the formal stable subspace 
$E^s(x)\subset\Cal T_x\bold M$ of $x$ as

$$
E^s(x)=\left\{(\delta q,\,\delta v)\in\Cal T_x\bold M\big|\; 
\delta v=-B(x)[\delta q]\right\},
$$
where the symmetric, positive semi-definite operator $B(x)$ (acting on the tangent space
of $\bold Q$ at the footpoint $q$, where $x=(q,v)$) is defined by the
continued fraction expansion introduced by Sinai in [Sin(1979)], see also
[Ch(1982)] or (2.4) in [K-S-Sz(1990)-I]. It is a well known fact that
$E^s(x)$ is the tangent space of the local stable manifold $\gamma^s(x)$, if
the latter object exists.

\medskip

For any phase point $x\in\bold M\setminus\partial\bold M$ with a non-singular
backward orbit $S^{(-\infty,0)}x$ (and with at least one collision on it) the
unstable space $E^u(x)$ of $x$ is defined as $-E^s(-x)$, where 
$-x=(q,\,-v)$ for $x=(q,\,v)$.

\medskip

All tangent vectors $(\delta q,\delta v)$ considered in Appendix II are not
just arbitrary tangent vectors: we restrict the exposition to the "orthogonal
section'', i. e. to the one-codimensional subspaces of the tangent space
consisting of all tangent vectors $(\delta q,\delta v)\in\Cal T_x\bold M$ for
which $\delta q$ is orthogonal to the velocity component $v$ of the phase
point $x=(q,\, v)$.

\bigskip

\subheading{Theorem} For any phase point
$x_0\in\bold M\setminus\partial\bold M$ with a non-singular
forward orbit $S^{(0,\infty)}x_0$ and with infinitely many consecutive,
connected collision graphs on $S^{(0,\infty)}x_0$, and for any
number $L>0$ one can find a time $t>0$ and a non-zero tangent vector
$(\delta q_0,\delta v_0)\in E^s(x_{0})$ with

$$
\frac{||(\delta q_t,\,\delta v_t)||}{||(\delta q_0,\,\delta v_0)||}<L^{-1},
$$
where $(\delta q_t,\delta v_t)=DS^t(\delta q_0,\delta v_0)\in E^s(x_t)$,
$x_t=S^tx_0$.

\medskip

\subheading{Proof}
Taking time-reversal, we would like to get useful lower estimates for the
expansion of a tangent vector $(\delta q_0,\delta v_0)\in\Cal T_{x_0}\bold M$
with positive infinitesimal Lyapunov function $Q(\delta q_0,\delta
v_0)=\langle\delta q_0,\delta v_0\rangle$. The expression $\langle\delta
q_0,\delta v_0\rangle$ is the scalar product in $\Bbb R^d$ defined via the
mass (or kinetic energy) metric, see \S2. It is also called the infinitesimal
Lyapunov function associated with the tangent vector $(\delta q_0,\delta
v_0)$, see [K-B(1994)], or part A.4 of the Appendix in [Ch(1994)], or \S7 of
[Sim(2003)]. For a detailed exposition of the relationship between the
quadratic form $Q(\,.\,)$, the relevant symplectic geometry of the Hamiltonian
system and the dynamics, please also see [L-W(1995)].

\medskip

\subheading{Note} The idea to use indefinite quadratic forms to study Anosov
systems belongs to Lewowicz, and in the case of non-uniformly hyperbolic
systems to Wojtkowski. In particular, the quadratic form $Q$ was introduced
into the subject of semi-dispersive billiards by Wojtkowski, implicitly in 
[W(1985)], explicitly in [W(1988)]. The symplectic formulation of the form
$Q$ appeared first in [W(1990)].

These ideas have been explored in detail and further developed by
N. I. Chernov and myself in recent personal communications, so that we
obtained at least linear (but uniform!) expansion rates for submanifolds
with negative infinitesimal Lyapunov forms for their normal vector. These
results are presented in our recent joint paper [Ch-Sim(2006)]. Also, closely
related to the above said, the following ideas (to estimate the expansion
rates of tangent vectors from below) are derived from the thoughts being
published in [Ch-Sim(2006)].

\medskip

Denote by $(\delta q_t,\delta v_t)=(DS^t)(\delta q_0,\delta v_0)$
the image of the tangent vector $(\delta q_0,\delta v_0)$ under the
linearization $DS^t$ of the map $S^t$. (We assume that the base phase
point $x_0$ --- for which $(\delta q_0,\delta v_0)\in\Cal T_{x_0}\bold M$ ---
has a non-singular forward orbit.) The time-evolution 
$(\delta q_{t_1},\delta v_{t_1})\mapsto(\delta q_{t_2},\delta v_{t_2})$
($t_1<t_2$) on a collision free segment $S^{[t_1,t_2]}x_0$ is described
by the equations

$$
\aligned
\delta v_{t_2}&=\delta v_{t_1}, \\
\delta q_{t_2}&=\delta q_{t_1}+(t_2-t_1)\delta v_{t_1}.
\endaligned
\tag A.1
$$
Correspondingly, the change 
$Q(\delta q_{t_1},\delta v_{t_1})\mapsto Q(\delta q_{t_2},\delta v_{t_2})$
in the infinitesimal Lyapunov function $Q(\,.\,)$ on the collision free orbit
segment $S^{[t_1,t_2]}x_0$ is

$$
Q(\delta q_{t_2},\delta v_{t_2})=Q(\delta q_{t_1},\delta v_{t_1})
+(t_2-t_1)||\delta v_{t_1}||^2,
\tag A.2
$$
thus $Q(\,.\,)$ steadily increases between collisions. 

The passage
$(\delta q_t^-,\delta v_t^-)\mapsto(\delta q_t^+,\delta v_t^+)$ through a
reflection (i. e. when $x_t=S^tx_0\in\partial\bold M$) is given by Lemma 2 of
[Sin(1979)] or formula (2) in \S3 of [S-Ch(1987)]:

$$
\aligned
\delta q_t^+&=R\delta q_t^-, \\
\delta v_t^+&=R\delta v_t^-+2\cos\phi RV^*KV\delta q_t^-,
\endaligned
\tag A.3
$$
where the operator $R:\, \Cal T\bold Q\to\Cal T\bold Q$ is the orthogonal 
reflection (with respect to the mass metric) across the tangent hyperplane
$\Cal T_{q_t}\partial\bold Q$ of the boundary $\partial\bold Q$ at the
configuration component $q_t$ of $x_t=(q_t,v_t^\pm)$, 
$V:\, (v_t^-)^\perp\to\Cal T_{q_t}\partial\bold Q$ is the $v_t^-$-parallel
projection of the orthocomplement hyperplane $(v_t^-)^\perp$ onto
$\Cal T_{q_t}\partial\bold Q$, 
$V^*:\,\Cal T_{q_t}\partial\bold Q\to(v_t^-)^\perp$ is the adjoint of $V$
(i. e. the $n(q_t)$-parallel projection of $\Cal T_{q_t}\partial\bold Q$
onto $(v_t^-)^\perp$, where $n(q_t)$ is the inner unit normal vector of 
$\partial\bold Q$ at $q_t\in\partial\bold Q$), 
$K:\, \Cal T_{q_t}\partial\bold Q\to\Cal T_{q_t}\partial\bold Q$ is the second
fundamental form of the boundary $\partial\bold Q$ at $q_t$ (with respect to 
the field $n(q)$ of inner unit normal vectors of $\partial\bold Q$) and,
finally, $\cos\phi=\langle n(q_t),\, v_t^+\rangle>0$ is the cosine of the 
angle $\phi$ ($0\le\phi<\pi/2$) subtended by $v_t^+$ and $n(q_t)$. Regarding
formulas (A.3), please see the last displayed formula in \S1 of 
[S-Ch(1987)] or (i)--(ii) in Proposition 2.3 of [K-S-Sz(1990)-I]. The
instanteneous change in the infinitesimal Lyapunov function
$Q(\delta q_t,\delta v_t)$ caused by the reflection at time $t$ is easily
derived from (A.3):

$$
\aligned
Q(\delta q_t^+,\delta v_t^+)&=Q(\delta q_t^-,\delta v_t^-)+2\cos\phi\langle
V\delta q_t^-,\, KV\delta q_t^-\rangle \\
&\ge Q(\delta q_t^-,\delta v_t^-).
\endaligned
\tag A.4
$$ 
In the last inequality we used the fact that the operator $K$ is positive
semi-definite, i. e. the billiard is semi-dispersive.

We are primarily interested in getting useful lower estimates for the expansion
rate $||\delta q_t||/||\delta q_0||$. The needed result is

\medskip

\subheading{Proposition A.5} Use all the notations above, and assume that

$$
\langle\delta q_0,\, \delta v_0\rangle/||\delta q_0||^2\ge c_0>0.
$$
We claim that $||\delta q_t||/||\delta q_0||\ge 1+c_0t$ for all $t\ge0$.

\medskip

\subheading{Proof} Clearly, the function $||\delta q_t||$ of $t$ is continuous
for all $t\ge0$ and continuously differentiable between collisions. According 
to (A.1), $\frac{d}{dt}\delta q_t=\delta v_t$, so

$$
\frac{d}{dt}||\delta q_t||^2=2\langle\delta q_t, \delta v_t\rangle.
\tag A.6
$$

Observe that not only the positive valued function 
$Q(\delta q_t,\delta v_t)=\langle\delta q_t,\delta v_t\rangle$ is 
nondecreasing in $t$ by (A.2) and (A.4), but the quantity 
$\langle\delta q_t,\delta v_t\rangle/||\delta q_t||$ is nondecreasing in $t$,
as well. The reason is that
$\langle\delta q_t,\delta v_t\rangle/||\delta q_t||=||\delta v_t||\cos\alpha_t$
($\alpha_t$ being the acute angle subtended by $\delta q_t$ and $\delta v_t$),
and between collisions the quantity $||\delta v_t||$ is unchanged, while the
acute angle $\alpha_t$ decreases, according to the time-evolution equations
(A.1). Finally, we should keep in mind that at a collision the norm
$||\delta q_t||$ does not change, while $\langle\delta q_t,\delta v_t\rangle$
cannot decrease, see (A.4). Thus we obtain the inequalities

$$
\langle\delta q_t,\delta v_t\rangle/||\delta q_t||\ge
\langle\delta q_0,\delta v_0\rangle/||\delta q_0||\ge c_0||\delta q_0||,
$$
so

$$
\frac{d}{dt}||\delta q_t||^2=2||\delta q_t||\frac{d}{dt}||\delta q_t||
=2\langle\delta q_t,\delta v_t\rangle
\ge 2c_0||\delta q_0||\cdot||\delta q_t||
$$
by (A.6) and the previous inequality. This means that 
$\frac{d}{dt}||\delta q_t||\ge c_0||\delta q_0||$, so
$||\delta q_t||\ge||\delta q_0||(1+c_0t)$, proving the proposition. \qed

\medskip

Next we need an effective lower estimate $c_0$ for the curvature
$\langle\delta q_0,\, \delta v_0\rangle/||\delta q_0||^2$ of the trajectory
bundle:

\medskip

\subheading{Lemma A.7} Assume that the perturbation 
$(\delta q_0^-,\,\delta v_0^-)\in\Cal T_{x_0}\bold M$ (as in Proposition A.5)
is being performed at time zero right before a collision, say,
$\sigma_0=(1,\,2)$ taking place at that time. Select the tangent vector
$(\delta q_0^-,\,\delta v_0^-)$ in such a specific way that $\delta
q_0^-=(m_2w,-m_1w,0,0,\dots,0)$ with a nonzero vector $w\in\Bbb R^\nu$,
$\langle w,v_1^--v_2^-\rangle=0$. This scalar product equation is exactly the
condition that guarantees that $\delta q_0^-$ be orthogonal to the velocity
component $v^-=(v_1^-,v_2^-,\dots,v_N^-)$ of $x_0=(q,v^-)$. The next, though
crucial requirement is that $w$ should be selected from the two-dimensional
plane spanned by $v_1^--v_2^-$ and $q_1-q_2$ (with $||q_1-q_2||=2r$) in $\Bbb
R^\nu$. If $v_1^- -v_2^-$ and $q_1-q_2$ are parallel, then we do not impose
this condition, for in that case there is no ``astigmatism''.
The purpose of this condition is to avoid the unwanted phenomenon of
``astigmatism'' in our billiard system, discovered first by Bunimovich and
Rehacek in [B-R(1997)] and [B-R(1998)]. Later on the phenomenon of astigmatism
gathered further prominence in the paper [B-Ch-Sz-T(2002)] as the main driving
mechanism behind the wild non-differentiability of the singularity manifolds
(at their boundaries) in semi-dispersive billiard systems with a configuration
space dimension bigger than $2$.
Finally, the last requirement is that the velocity component $\delta
v_0^-$ (right before the collision $(1,2)$) is chosen in such a way that the
tangent vector $(\delta q_0^-,\,\delta v_0^-)$ belongs to the unstable space
$E^u(x_0)$ of $x_0$. This can be done, indeed, by taking $\delta
v_0^-=B^u(x_0)[\delta q_0^-]$, where $B^u(x_0)$ is the curvature operator of the
unstable manifold of $x_0$ at $x_0$, right before the collision $(1,2)$ taking
place at time zero.

We claim that

$$
\frac{\langle\delta q_0^+,\delta v_0^+\rangle}{||\delta q_0||^2}\ge
\frac{||v_1-v_2||}{r\cos\phi_0}\ge\frac{||v_1-v_2||}{r}
\tag A.8
$$
for the post-collision tangent vector $(\delta q_0^+,\delta v_0^+)$, where
$\phi_0$ is the acute angle subtended by $v_1^+-v_2^+$ and the outer normal
vector of the sphere $\left\{y\in\Bbb R^\nu\big|\; ||y||=2r\right\}$ at the 
point $y=q_1-q_2$. Note that in (A.8) there is no need to use $+$ or $-$ in 
$||\delta q_0||^2$ or $||v_1-v_2||$, for $||\delta q_0^-||=||\delta q_0^+||$,
$||v_1^--v_2^-||=||v_1^+-v_2^+||$.

\medskip

\subheading{Proof} The proof of the first inequality in (A.8) is a simple, elementary
geometric argument in the plane spanned by $v_1^--v_2^-$ and $q_1-q_2$, so we
omit it. We only note that the outgoing relative velocity $v_1^+-v_2^+$ is
obtained from the pre-collision relative velocity $v_1^--v_2^-$ by reflecting
the latter one across the tangent hyperplane of the sphere 
$\left\{y\in\Bbb R^\nu\big|\; ||y||=2r\right\}$ at the point $y=q_1-q_2$.
It is a useful advice, though, to prove the first inequality of (A.8) in the
case $\delta v_0^-=0$ first (this is an elementary geometry exercise), then
observe that this inequality can only be further improved when we replace 
$\delta v_0^-=0$ with $\delta v_0^-=B^u(x_0)[\delta q_0^-]$. \qed

\medskip

The previous lemma shows that, in order to get useful lower estimates for the
``curvature'' $\langle\delta q,\delta v\rangle/||\delta q||^2$ of the
trajectory bundle, it is necessary
(and sufficient) to find collisions $\sigma=(i,j)$ on the orbit of a given 
point $x_0\in\bold M$ with a ``relatively big'' value of $||v_i-v_j||$.
Finding such collisions will be based upon the following result:

\medskip

\subheading{Proposition A.9} Consider orbit segments $S^{[0,T]}x_0$ of 
$N$-ball systems with masses $m_1,m_2,\dots,m_N$ in $\Bbb T^\nu$ (or in
$\Bbb R^\nu$) and with collision sequences 
$\Sigma=(\sigma_1,\sigma_2,\dots,\sigma_n)$ corresponding to connected 
collision graphs. (Now the kinetic energy is not necessarily normalized, and 
the total momentum $\sum_{i=1}^N m_iv_i$ may be different from zero.) We claim
that there exists a positive-valued function $f(a;m_1,m_2,\dots,m_N)$ 
($a>0$, $f$ is independent of the orbit segments $S^{[0,T]}x_0$) with the
following two properties:

\medskip

(1) If $||v_i(t_l)-v_j(t_l)||\le a$ for all collisions $\sigma_l=(i,j)$
($1\le l\le n$, $t_l$ is the time of $\sigma_l$) for some trajectory segment 
$S^{[0,T]}x_0$ with a symbolic collision sequence
$\Sigma=(\sigma_1,\sigma_2,\dots,\sigma_n)$ corresponding to a connected 
collision graph, then the norm $||v_{i'}(t)-v_{j'}(t)||$ of any relative
velocity at any time $t\in\Bbb R$ is at most $f(a;m_1,\dots,m_N)$;

(2) $\lim_{a\to 0} f(a;m_1,\dots,m_N)=0$ for any $(m_1,\dots,m_N)$.

\medskip

\subheading{Proof} We begin with

\medskip

\subheading{Lemma A.10} Consider an $N$-ball system with masses $m_1,\dots,m_N$
(an $(m_1,\dots\allowmathbreak,m_N)$-system, for short) in $\Bbb T^\nu$ 
(or in $\Bbb R^\nu$).
Assume that the inequalities $||v_i(0)-v_j(0)||\le a$ hold true 
($1\le i<j\le N$) for all relative velocities at time zero. We claim that 

$$
||v_i(t)-v_j(t)||\le 2a\sqrt{\frac{M}{m}}
\tag A.11
$$
for any pair $(i,j)$ and any time $t\in\Bbb R$, where $M=\sum_{i=1}^N m_i$ and 

$$
m=\min\left\{m_i|\; 1\le i\le N\right\}.
$$

\medskip

\subheading{Note} The estimate (A.11) is not optimal, however, it
will be sufficient for our purposes.

\medskip

\subheading{Proof} The assumed inequalities directly imply (by using a simple
convexity argument) that $||v'_i(0)||\le a$ ($1\le i\le N$) for the velocities
$v'_i(0)$ measured at time zero in the baricentric reference
system. Therefore, for the total kinetic energy $E_0$ (measured in the
baricentric system) we get the upper estimation $E_0\le \frac{1}{2}Ma^2$, and
this inequality remains true at any time $t$.  This means that all the
inequalities $||v'_i(t)||^2\le\frac{M}{m_i}a^2$ hold true for the baricentric
velocities $v'_i(t)$ at any time $t$, so

$$
||v'_i(t)-v'_j(t)||\le a\sqrt{M}\left(m_i^{-1/2}+m_j^{-1/2}\right)\le
2a\sqrt{\frac{M}{m}},
$$
thus the inequalities

$$
||v_i(t)-v_j(t)||\le 2a\sqrt{\frac{M}{m}}
$$
hold true, as well. \qed

\medskip

\subheading{Proof of the proposition by induction on the number $N$}

\medskip

For $N=1$ we can take $f(a;m_1)=0$, and for $N=2$ the function
$f(a;m_1,m_2)=a$ is obviously a good choice for $f$. Let $N\ge 3$, and assume 
that the orbit segment $S^{[0,T]}x_0$ of an $(m_1,\dots,m_N)$-system fulfills
the conditions of the proposition. Let $\sigma_k=(i,j)$ be the collision in the
symbolic sequence $\Sigma_n=(\sigma_1,\dots,\sigma_n)$ of $S^{[0,T]}x_0$ with
the property that the collision graph of $\Sigma_k=(\sigma_1,\dots,\sigma_k)$
is connected, while the collision graph of 
$\Sigma_{k-1}=(\sigma_1,\dots,\sigma_{k-1})$ is still disconnected. Denote the
two connected components (as vertex sets) of $\Sigma_{k-1}$ by $C_1$ and
$C_2$, so that $i\in C_1$, $j\in C_2$, $C_1\cup C_2=\{1,2,\dots,N\}$, and
$C_1\cap C_2=\emptyset$. By the induction hypothesis and the condition of the
proposition, the norm of any relative velocity
$v_{i'}(t_k-0)-v_{j'}(t_k-0)$ right before the collision $\sigma_k$ (taking 
place at time $t_k$) is at most 
$a+f(a;\,\overline{C}_1)+f(a;\,\overline{C}_2)$, where $\overline{C}_l$ stands
for the collection of the masses of all particles in the component $C_l$,
$l=1,\,2$. Let $g(a;m_1,\dots,m_N)$ be the maximum of all possible sums

$$
a+f(a;\,\overline{D}_1)+f(a;\,\overline{D}_2),
$$
taken for all two-class partitions $(D_1,D_2)$ of the vertex set
$\{1,2,\dots,N\}$. According to the previous lemma, the function

$$
f(a;m_1,\dots,m_N):=2\sqrt{\frac{M}{m}}g(a;m_1,\dots,m_N)
$$
fulfills both requirements (1) and (2) of the proposition. \qed

\medskip

\subheading{Corollary A.12} Consider the original $(m_1,\dots,m_N)$-system with
the standard normalizations $\sum_{i=1}^N m_iv_i=0$, 
$\frac{1}{2}\sum_{i=1}^N m_i||v_i||^2=\frac{1}{2}$. We claim that there exists
a threshold $G=G(m_1,\dots,m_N)>0$ (depending only on $N$, $m_1,\dots,m_N$)
with the following property:

In any orbit segment $S^{[0,T]}x_0$ of the $(m_1,\dots,m_N)$-system with the
standard normalizations and with a connected collision graph, one can always
find a collision $\sigma=(i,j)$, taking place at time $t$, so that
$||v_i(t)-v_j(t)||\ge G(m_1,\dots,m_N)$.

\medskip

\subheading{Proof} Indeed, we choose $G=G(m_1,\dots,m_N)>0$ so small that
$f(G;m_1,\dots,m_N)\allowmathbreak<M^{-1/2}$. 
Assume, contrary to A.12, that the norm of any
relative velocity $v_i-v_j$ of any collision of $S^{[0,T]}x_0$ is less than the
above selected value of $G$. By the proposition, we have the inequalities
$||v_i(0)-v_j(0)||\le f(G;m_1,\dots,m_N)$ at time zero. The normalization
$\sum_{i=1}^N m_iv_i(0)=0$, with a simple convexity argument, implies that
$||v_i(0)||\le f(G;m_1,\dots,m_N)$ for all $i$, $1\le i\le N$, so the total
kinetic energy is at most
$\frac{1}{2}M\left[f(G;m_1,\dots,m_N)\right]^2<\frac{1}{2}$, a contradiction.
\qed

\medskip

\subheading{Corollary A.13} For any phase point $x_0$ with a non-singular
backward trajectory $S^{(-\infty,0)}x_0$ and with infinitely many consecutive,
connected collision graphs on $S^{(-\infty,0)}x_0$, and for any number $L>0$
one can always find a time $-t<0$ and a non-zero tangent vector $(\delta
q_0,\delta v_0)\in E^u(x_{-t})$ ($x_{-t}=S^{-t}x_0$) with $||\delta
q_t||/||\delta q_0||>L$, where $(\delta q_t,\delta v_t)=DS^t(\delta q_0,\delta
v_0)\in E^u(x_0)$.

\medskip

\subheading{Proof} Indeed, select a number $t>0$ so big that
$1+\frac{t}{r}G(m_1,\dots,m_N)>L$ and $-t$ is the time of a collision (on the
orbit of $x_0$) with the relative velocity $v^-_i(-t)-v^-_j(-t)$, for which
$||v^-_i(-t)-v^-_j(-t)||\ge G(m_1,\dots,m_N)$. Indeed, this can be done,
thanks to A.12 and the assumed abundance of connected collision graphs.
By Lemma A.7 we can choose a
non-zero tangent vector $(\delta q^-_0,\,\delta v^-_0)\in E^u(x_{-t})$ right
before the collision at time $-t$ in such a way that the lower estimate

$$
\frac{\langle\delta q^+_0,\,\delta v^+_0\rangle}{||\delta q^+_0||^2}\ge
\frac{1}{r}G(m_1,\dots,m_N)
$$
holds true for the ``curvature''
$\langle\delta q^+_0,\,\delta v^+_0\rangle/||\delta q^+_0||^2$ associated
with the post-collision tangent vector $(\delta q^+_0,\,\delta v^+_0)$.
According to Proposition A.5, we have the lower estimate

$$
\frac{||\delta q_t||}{||\delta q_0||}\ge 1+\frac{t}{r}G(m_1,\dots,m_N)>L
$$
for the $\delta q$-expansion rate between 
$(\delta q^+_0,\, \delta v^+_0)$ and
$(\delta q_t,\,\delta v_t)=DS^t(\delta q^+_0,\, \delta v^+_0)$. \qed

\medskip

We remind the reader that, according to the main result of [B-F-K(1998)],
there exists a number $\epsilon_0=\epsilon_0(m_1,\dots,m_N;\,r;\,\nu)>0$
and a large threshold $N_0=N_0(m_1,\dots,m_N;\,r;\,\nu)\in\Bbb N$ such that in
the $(m_1,\dots,m_N;\,r;\,\nu)$-billiard flow amongst any $N_0$ consecutive
collisions one can always find two neighboring ones separated in time by at
least $\epsilon_0$. Thus, for a phase point $x_{-t}$ at least
$\epsilon_0/2$-away from collisions, the norms $||\delta q_0||$ and
$\sqrt{||\delta q_0||^2+||\delta v_0||^2}$ are equivalent for all vectors
$(\delta q_0,\,\delta v_0)\in E^u(x_{-t})$, hence, from the previous corollary
we immediately get

\medskip

\subheading{Corollary A.14}
For any phase point $x_0\in\bold M\setminus\partial\bold M$ with a
non-singular backward trajectory $S^{(-\infty,0)}x_0$ and with infinitely many
consecutive, connected collision graphs on $S^{(-\infty,0)}x_0$, and for any
number $L>0$ one can always find a time $-t<0$ and a non-zero tangent vector
$(\delta q_0,\delta v_0)\in E^u(x_{-t})$ ($x_{-t}=S^{-t}x_0$) with

$$
\frac{||(\delta q_t,\,\delta v_t)||}{||(\delta q_0,\,\delta v_0)||}>L,
$$
where 
$(\delta q_t,\delta v_t)=DS^t(\delta q_0,\delta v_0)\in E^u(x_0)$. \qed

\medskip

The time-reversal dual of the previous result is immediately obtained by
replacing the phase point $x_0=(q_0,\,v_0)$ with $-x_0=(q_0,\,-v_0)$, the
backward orbit with the forward orbit, and the unstable vectors with the
stable ones. This dual result is exactly our theorem, formulated at the
beginning of Appendix II.

\bigskip

\subheading{Acknowledgement} The author expresses his sincere gratitude
to N. I. Chernov for his numerous, very valuable questions, remarks, and
suggestions.

\bye